\documentclass[12pt]{article}

\usepackage{amssymb,amsmath}

\usepackage{color}
\usepackage{graphicx}

\def\a{{\bf a}}

\def\f{{\bf f}}

\def\i{{\bf i}}
\def\j{{\bf j}}

\def\u{{\bf u}}
\def\v{{\bf v}}

\def\x{{\bf x}}

\def\E{{\cal E}}

\def\G{{\cal G}}

\def\M{{\cal M}}

\def\Abb{{\mathbb A}}
\def\Bbb{{\mathbb B}}

\def\Nbb{{\mathbb N}}

\def\R{{\mathbb R}}
\def\Sm{{\mathbb S}}

\def\d{\delta}

\def\GA{\Gamma}

\def\t{\tau}
\def\th{\theta}

\def\bmu{{\boldsymbol \mu}}
\def\bnu{{\boldsymbol \nu}}
\def\bth{{\boldsymbol \theta}}
\def\bphi{{\boldsymbol \phi}}
\def\bpsi{{\boldsymbol \psi}}
\def\bpi{{\boldsymbol \pi}}

\def\bzeta{{\boldsymbol \zeta}}

\def\Pt{\tilde{P}}

\def\kai{k \ap \infty}
\def\lai{l \ap \infty}

\def\ap{\rightarrow}

\def\seq{\subseteq}

\def\bz{{\bf 0}}

\def\fa{\; \forall}

\def\st{\mbox{ s.t. }}

\def\nm{\Vert}

\renewcommand{\and}{\mbox{$\wedge$}}

\newcommand{\bc}{\begin{center}}
\newcommand{\ec}{\end{center}}
\newcommand{\be}{\begin{equation}}
\newcommand{\ee}{\end{equation}}
\newcommand{\bd}{\begin{displaymath}}
\newcommand{\ed}{\end{displaymath}}
\newcommand{\ba}{\begin{array}}
\newcommand{\ea}{\end{array}}
\newcommand{\ben}{\begin{enumerate}}
\newcommand{\een}{\end{enumerate}}
\newcommand{\bit}{\begin{itemize}}
\newcommand{\eit}{\end{itemize}}
\newcommand{\beq}{\begin{eqnarray}}
\newcommand{\eeq}{\end{eqnarray}}
\newcommand{\btab}{\begin{tabular}}
\newcommand{\etab}{\end{tabular}}
\newcommand{\bfig}{\begin{figure}}
\newcommand{\efig}{\end{figure}}

\def\halmos{\hfill{\blacksquare}}

\def\m{{\bf m}}

\def\bmubar{\bar{\bmu}}
\def\bnubar{\bar{\bnu}}
\def\bthbar{\bar{\bth}}
\def\bzetabar{\bar{\bzeta}}

\def\mubar{\bar{\mu}}
\def\nubar{\bar{\nu}}
\def\thbar{\bar{\th}}
\def\zetabar{\bar{\zeta}}

\def\lbar{\bar{l}}

\def\bmuhat{\hat{\bmu}}

\def\GA{\Gamma}
\def\MA{\M(\Abb)}
\def\McA{\M_s(\Abb^2)}
\def\McAs{\M_s(\Abb^{s+1})}

\newtheorem{corollary}{Corollary}{\bf}{\it}
\newtheorem{definition}{Definition}{\bf}{\it}
{\bf}{\rm}
\newtheorem{lemma}{Lemma}{\bf}{\it}
\newtheorem{theorem}{Theorem}{\bf}{\it}
{\bf}{\it}
{\bf}{\it}
{\bf}{\rm}

\begin{document}


\title{An Elementary Derivation of the\\
Large Deviation Rate Function for\\
Finite State Markov Chains}
\author{M.\ Vidyasagar
\thanks{The author is Cecil \& Ida Green Chair,
Erik Jonsson School of Engineering \& Computer Science,
University of Texas at Dallas, 800 W.\ Campbell Road, Richardson, TX
75080, USA; email: M.Vidyasagar@utdallas.edu.
This research was supported by the National Science Foundation of the USA
under Award \# 1001643 and the Cecil \& Ida Green Endowment to the University
of Texas at Dallas.}}
\date{}
\maketitle

\begin{abstract}
Large deviation theory is a branch of probability theory that
is devoted to a study of the ``rate'' at which empirical estimates
of various quantities converge to their true values.
The object of study in this paper is the rate at which estimates
of the doublet frequencies of a Markov chain over a finite alphabet
converge to their true values.
In case the Markov process is actually an i.i.d.\ process,
the rate function turns out to be the relative entropy (or Kullback-Leibler
divergence) between the true and the estimated probability vectors.
This result is a special case of a very general result
known as Sanov's theorem and dates back to 1957.
Moreover, since the introduction of the ``method of types'' by
Csisz\'{a}r and his co-workers during the 1980s, the proof of
this version of Sanov's theorem has been ``elementary,''
using some combinatorial arguments.
However, when the i.i.d.\ process is replaced by a Markov process,
the available proofs are far more complex.
The main objective of this paper is therefore to present
a first-principles derivation of the LDP for finite state Markov chains,
using only simple combinatorial arguments (e.g.\ the method of types),
thus gathering in one place various arguments and estimates that
are scattered over the literature.
\end{abstract}

\section{Introduction}\label{sec:1}

Large deviation theory is a branch of probability theory that
is devoted to a study of the ``rate'' at which empirical estimates
of various quantities converge to their true values.
Unlike in statistical learning theory (see e.g.\
\cite{Vapnik98,MV-Book97,MV-Book02}), where the emphasis is on
deriving ``finite-time'' estimates, large deviation theory is only
asymptotic and tells us what happens in the limit as the number
of samples approaches infinity.

As the name implies, the main objective of large deviation theory is
to estimate the likelihood that an empirical estimate of some
quantity exhibits a ``large deviation'' (i.e., differs significantly)
from its true value.
For instance (and this was one of the original motivations for studying
large deviation theory), suppose an insurance company both receives
life insurance premiums and also pays out as its policy-holders die.
How much money does it need to keep on hand in order to cover a possible
excess of claims over income?
Clearly, the company would like to keep as little reserve as possible,
since usually such reserve does not earn any return.
But if the claims in one particular period
exceed income plus reserve, the company may go bankrupt.
Therefore the company has an interest in estimating the level of reserve
that is sufficient to ensure that the probability of bankruptcy
does not exceed a prespecified threshold, say 1\%.
Now, based on  historical averages, the company may have available
the probability distribution of the mortality of its clients and thus
the possible claims, which may be called the ``true'' probability distribution.
However, during a particular period, it is possible for the actual mortality
and therefore the actual claims to be
substantially higher than the historical average or the ``true'' distribution.
Estimating the likelihood that such a ``large deviation'' occurs is
the objective of the theory.

The subject of large deviation theory
in one form or another dates back to the 1920s if not earlier.
Some of the earliest work (motivated by the actuarial problem
mentioned in the previous paragraph) was done by Cram\'{e}r
during the 1930s; see \cite{Cramer36,Cramer38}.
Cram\'{e}r's work is applicable to real-valued random variables,
but is based on the assumption that successive samples
are independent.
Much later, the work of G\"{a}rtner \cite{Gartner77} and Ellis
\cite{Ellis84} removed the assumption that samples are independent,
and thus extended the applicability of the results to much more
general settings, including the situation where the
samples come from a Markov process.
For i.i.d.\ samples, one of the main results, known as Sanov's
theorem \cite{Sanov57}, dates to 1957.
For random variables assumming values in a finite set,
the theory was considerably simplified by Csisz\'{a}r and his co-workers
in the late 1970s and early 1980s, via the introduction of the so-called
``method of types.''
As a result, large deviation theory for i.i.d.\ (independent identically
distributed) processes over a finite alphabet is now well-understood,
and the proofs of the main results are easily accessible.

In contrast, currently available treatments of large deviation theory
for Markov chains over a finite alphabet are far more involved, 
even in widely used research monographs.
For instance, in \cite{DZ98}, the rate function (to be defined later)
for i.i.d.\ processes over a finite alphabet is given on page 16,
while the rate function for Markov processes over a finite alphabet
is given on page 79, thus suggesting that more than
sixty pages of mathematics are required before one can make the
transition from i.i.d.\ processes to Markov chains.
In \cite{Hollander00} the proof is more direct than in \cite{DZ98}, but
at a crucial point it invokes Varadhan's lemma, which is a very advanced
concept and is not really needed to study Markov chains over a finite alphabet.
Historically, large deviation theory for Markov processes is contained
in the work of Donsker and Varadhan \cite{DV75a,DV75b,DV75c,DV76,DV83}.
Extensions to still more general situations were carried out by G\"{a}rtner
\cite{Gartner77} and Ellis \cite{Ellis84}, leading to a result known
as the G\"{a}rtner-Ellis theorem; see for example \cite[p. 34]{DZ98}.
Indeed, in \cite{DZ98} large deviation theory for Markov chains
is derived using the G\"{a}rtner-Ellis theorem as the starting point.
Consequently, the proof is rather round-about and not so simple as it could be.

The complexity of the treatments of Markov chains in various books is
a little surprising, because simpler proofs are indeed available
in the literature, though they are scattered.
For instance, a very old paper by Natarajan \cite{Natarajan85}
explicitly extends the method of types to Markov chains, and makes use
of some bounds on the size of type classes derived earlier by Davisson
et al. \cite{DLS81}.
However, the paper of Natarajan does not seem to be so well-known;
indeed neither \cite{Natarajan85} nor \cite{DLS81} is cited in \cite{DZ98},
which has a very extensive bibliography.
There is another paper by Csisz\'{a}r et al. \cite{CCC87} that
discusses the large deviation property for Markov chains in a slightly
roundabout fashion, whereas the paper \cite{Natarajan85} has the derivation of
the rate function as its principal objective.

This paper is motivated by the present author's belief that 
it should be possible to prove the basic results of large deviation theory
for Markov chains over a finite alphabet
almost as easily as for i.i.d.\ processes.
Towards this end, several simplifications that are currently scattered
in the literature are collected in one place and streamlined.
As such there are no ``new results'' in the paper.
Rather, the hope is that by making the theory accessible without invoking
advanced concepts, the paper would encourage the research community to
make greater use of large deviation theory in the context of Markov chains.

One of the advantages of the first-principles derivation given here is that
the extension to multi-step Markov chains is straight-forward.
Indeed, through some elementary computations,
the rate function for Markov chains as derived 
in early works \cite{DLS81,Natarajan85} can be shown to be just the
``conditional relative entropy.''
This interpretation was present in early works \cite{DLS81,Natarajan85} but
is not emphasized in today's literature.
In the opinion of the author, the interpretation of the rate function as the
conditional relative entropy is very natural, because it can be readily
extended to multi-step Markov chains, and can also be extended with some
technical assumptions to stationary stochastic processes that have
a fading memory, as is attempted in \cite{CG05}.

\section{Preliminaries}\label{sec:2}

In this section we briefly review of some concepts from information
theory that are needed in subsequent sections.
This review also serves to fix notation.
A reader who is encountering these concepts for the first time should
consult any standard text, e.g.\ \cite{CT06}, for further details.

Throughout we will study Markov chains assuming values in a finite set
$\Abb$ of cardinality $n$.
Strictly speaking, we should denote the elements of $\Abb$ as $\{ a_1 ,
\ldots , a_n \}$, but in the interests of simplicity we will write
$\Abb = \{ 1 , \ldots , n \}$.
The reader should remember however that these are just ``symbols'' or
``labels,'' and do not correspond to the integers from $1$ to $n$.

\subsection{Stationary Distributions and an Alternate Description of
Markov Chains}\label{ssec:21}

Let $\Sm_n$ denote the $n$-dimensional simplex, i.e.
\bd
\Sm_n := \{ \v \in [0,1]^n : \sum_{i=1}^n v_i = 1 \} .
\ed
Thus if $\Abb$ is a set of cardinality $n$, then the set of all
probability measures on $\Abb$, denoted by $\M(\Abb)$, can be
identified with $\Sm_n$.
Similarly, if $X$ is a random variable assuming values in a finite set
$\Abb = \{ 1 , \ldots , n \}$, the associated distribution
\bd
\bmu := [ \Pr \{ X = i \} , i = 1 , \ldots , n ]
\ed
belongs to $\Sm_n$.
Throughout the paper, we use bold face Greek letters to denote probability
distributions, that is, vectors in a simplex of appropriate dimension.

Next, suppose $\{ X_t \}$ is a stationary Markov process assuming values
in $\Abb$.
The customary way to describe such a Markov chain is in terms of
two entities:
(i) The {\bf stationary distribution} $\bpi \in \MA$, where
\bd
\pi_i := \Pr \{ X_t = i \} , i \in \Abb ,
\ed
and (ii) The {\bf state transition matrix} $A \in [0,1]^{n \times n}$, where
\bd
a_{ij} := \Pr \{ X_{t+1} = j | X_t = i \} , \fa i, j \in \Abb .
\ed
Since the process is stationary, all of the above quantities do not depend
on the time $t$.
Note that each row of $A$ is a probability distribution and belongs to
$\MA$, since the $i$-th row of $A$ is the conditional distribution of
$X_{t+1}$ given that $X_t = i$.
Thus $A$ is a stochastic matrix.
Note that it is in general not enough to specify the transition matrix $A$
alone, because in general there could be more than one stationary distribution
$\bpi \in \Sm_n$ that satisfies $\bpi = \bpi A$.

But there is an alternative description of a Markov chain that is more 
convenient for present purposes, namely the {\em vector of doublet frequencies}.
To motivate this alternate description, we first introduce the notion
of a stationary distribution on $\Abb^k$.\footnote{The present author
actually prefers the phrase ``consistent'' distribution, but unfortunately
``stationary'' is far more commonly used.}

Suppose $\{ X_t \}$ is a stationary stochastic process (not necessarily
Markov) assuming values in a finite set $\Abb$.
To make notation more compact, define
\bd
X_s^t := X_s X_{s+1} \ldots X_{t-1} X_t .
\ed
Clearly this notation makes sense only when $s \leq t$.
For each integer $k \geq 1$, and each $\i := (i_1 , \ldots , i_k) \in \Abb^k$,
let us define the $k$-tuple frequencies
\bd
\mu_\i := \Pr \{ X_{t+1}^{t+k} = \i \} .
\ed
Again, this probability does not depend on $t$ since the process is stationary.
Now note that, for each $k$-tuple $\i \in \Abb^k$, the events
\bd
\{ X_{t+1}^{t+k+1}  = (\i,1) \} , \ldots , \{ X_{t+1}^{t+k+1}  = (\i,n) \} 
\ed
are pairwise disjoint, and together generate the event
\bd
\{ X_{t+1}^{t+k} = \i \} .
\ed
Thus
\bd
\mu_\i = \sum_{j \in \Abb} \mu_{\i j} , \fa \i \in \Abb^k .
\ed
By entirely analogous reasoning, it also follows that
\bd
\mu_\i = \sum_{j \in \Abb} \mu_{j \i} , \fa \i \in \Abb^k .
\ed
This motivates the next definition.

\begin{definition}\label{def:1}
A distribution $\bnu \in \M(\Abb^2)$ is said to be {\bf stationary}
(or consistent) if
\be\label{eq:21}
\sum_{j \in \Abb} \nu_{ij} = \sum_{j \in \Abb} \nu_{ji} , \fa i \in \Abb .
\ee
For $k \geq 2$, a distribution $\bnu \in \M(\Abb^{k+1})$ is said to be
{\bf stationary} (or consistent) if
\bd
\sum_{j \in \Abb} \nu_{\i j} = \sum_{j \in \Abb} \nu_{j \i} ,
\fa \i \in \Abb^k,
\ed
and in addition, the resulting distribution $\bnubar$ on $\Abb^k$
defined by
\bd
\bnubar_{\i} := \sum_{j=1}^n \nu_{\i j} \fa \i \in \Abb^k
= \sum_{j=1}^n \nu_{j \i} \fa \i \in \Abb^k
\ed
is stationary.
Equivalently, a distribution $\bnu \in \M(\Abb^k)$ is said to be
{\bf stationary} (or consistent) if

\be\label{eq:22}
\sum_{ \j_1 \in \Abb^{l_1} } \sum_{ \j_2 \in \Abb^{l_2} }
\nu_{\j_1 \i \j_2} = \nu_\i ,
\ee
for all $\i \in \Abb^{k - l_1 - l_2}$, and all $l_1 \geq 0, l_2 \geq 0$ with
$l_1 + l_2 \leq k$.
\end{definition}
In the above definition, if $l_1 = 0$, then the corresponding summation
is removed, and similarly if $l_2 = 0$.

Note that some authors define stationarity through (\ref{eq:21}) alone
and don't impose the additional restriction that the resulting reduced
distribution must also be stationary.
This is clearly incorrect.
We use the symbol $\M_s(\Abb^k)$ to denote the set of all stationary
distributions on $\Abb^k$.
Clearly $\M_s(\Abb^k)$ is a subset of $\M(\Abb^k)$, the set of {\it all\/}
distributions on $\Abb^k$, which is the same as $ \Sm_{n^k}$.

Now let us return to Markov chains.
Suppose $\{ X_t \}$ is a stationary Markov chain assuming values in the
finite set $\Abb$.
Define the vector $\bmu \in \M(\Abb^2)$ by
\be\label{eq:22a}
\mu_{ij} = \Pr \{ X_t X_{t+1} = ij \} , \fa i, j \in \Abb .
\ee
Then, as per the above discussion, actually $\bmu \in \McA$.
The claim is that the doublet frequency vector $\bmu$ captures all
the relevant information about the Markov chain.
The stationary distribution of the Markov chain is given by
\be\label{eq:22b}
\mubar_i := \sum_{j \in \Abb} \mu_{ij} = \sum_{j \in \Abb} \mu_{ji} ,
\ee
while the state transition matrix $A$ is given by
\bd
a_{ij} = \frac{\mu_{ij}}{\mubar_i} .
\ed
Dividing by $\mubar_i$ can be justified by observing that if $\mubar_i
= 0$ for some index $i$, then the corresponding element $i$ is never
visited by the Markov chain, and can thus be dropped from the set $\Abb$.
With these definitions, it readily follows that $\bmubar$ is a row
eigenvector of $A$, because
\bd
( \bmubar A )_j = \sum_{i=1}^n \mubar_i a_{ij} 
= \sum_{i=1}^n \mu_{ij} = \mubar_j .
\ed
Note that the above reasoning breaks down if
$\bmu \in \M(\Abb^2)$ but $\bmu \not \in \McA$.

More generally, suppose $\{ X_t \}$ is an $s$-step Markov chain, so that
\bd
\Pr \{ X_t | X^{t-1}_0 \} = \Pr \{ X_t | X^{t-1}_{t-s} \}
\fa t \geq s .
\ed
In other words, whenever $t \geq s$, the conditional distribution of
$X_t$ given $X_0^{t-1}$ is the same as the conditional distribution of
$X_t$ given the previous $s$ states.
Then the process is completely characterized by its $(s+1)$-tuple frequencies
\be\label{eq:22c}
\mu_\i := \Pr \{ X_t^{t+s} = \i \} , \fa \i \in \Abb^{s+1} .
\ee
The probability distribution $\bmu$ is stationary and thus belongs to $\McAs$.
Now an $s$-step Markov chain assuming values in $\Abb$ can also be viewed
as a conventional (one-step) Markov chain over the state space $\Abb^s$.
Moreover, if the current state is $i \j$ where $i \in \Abb , \j \in 
\Abb^{s-1}$, then a transition is possible only to a state of the form
$\j k, k \in \Abb$.
Thus, even though the state transition matrix has dimension $n^s \times n^s$,
each row of the transition matrix can have at most $n$ nonzero elements.
In row $i \j$, the entry in column $\j k$ equals
\bd
\Pr \{ X_t = k | X_{t-s}^{t-1} = i \j \} = \frac{ \mu_{i \j k} }
{ \mubar_{i \j} },
\ed
where as before we define
\be\label{eq:22d}
\mubar_{{\bf k}} := \sum_{i \in \Abb} \mu_{i {\bf k}} 
= \sum_{i \in \Abb} \mu_{{\bf k} i } , \fa {\bf k} \in \Abb^s ,
\ee
while the entries in all other columns are zero.

\subsection{Entropy, Relative Entropy and Conditional Entropy}\label{ssec:22} 

In this subsection we review the definitions of various forms of entropy
for random variables.
This material lays the groundwork for the discussion of stationary
stochastic processes in the next subsection.

Given a distribution $\bnu \in \Sm_n$, its {\bf entropy} is defined as
\bd
H(\bnu) := \sum_{i=1}^n \nu_i \log(1/\nu_i) = - \sum_{i=1}^n \nu_i \log \nu_i ,
\ed
where $0 \log 0$ is taken as $0$.
If $\bmu, \bnu \in \Sm_n$, then their {\bf relative entropy} (also known
as the {\bf Kullback-Leibler divergence}) is defined as
\bd
D( \bnu \nm \bmu ) := \sum_{i=1}^n \nu_i \log \left( \frac{\nu_i}{\mu_i}
\right) .
\ed
Note that in order for $D( \bnu \nm \bmu )$ to be finite, it must be the
case that if $\mu_i = 0$ for some index $i$, then $\nu_i = 0$.
In such a case we say that $\bmu$ is ``dominated'' by $\bnu$ 
or that $\bnu$ ``dominates'' $\bmu$, and write
$\bmu \ll \bnu$ or $\bnu \gg \bmu$.
Note that the relative entropy can also be expressed as
\be\label{eq:23}
D( \bnu \nm \bmu ) = J( \bnu , \bmu ) - H( \bnu ),
\ee
where $J( \bnu , \bmu )$ is the ``loss function'' defined by
\be\label{eq:24}
J( \bnu , \bmu ) := \sum_{i=1}^n \nu_i \log ( 1/ \mu_i ) .
\ee
It is well-known that the relative entropy function $D( \cdot \nm \cdot )$
is jointly strictly convex in both arguments; see e.g.\ \cite{CT06}, Theorem
2.7.2.
As a result, for each fixed $\bmu$ the function $D( \bnu \nm \bmu )$ is a
strictly convex function of $\bnu$.

The next notion of the ``conditional entropy'' of one random variable
with respect to another.
Suppose $X,Y$ are random variables assuming values in finite sets
$\Abb = \{ 1 , \ldots , n\}$ and $\Bbb = \{ 1 , \ldots , m \}$ respectively.
Let $\bth \in \Sm_{nm}$ denote their joint distribution, that is
\bd
\th_{ij} = \Pr \{ X = i \& Y = j \} .
\ed
Then it is easy to see that
\bd
\Pr \{ X = i \} = \sum_{j=1}^m \Pr \{ X = i \& Y = j \} 
= \sum_{j=1}^m \th_{ij} ,
\ed
and analogously
\bd
\Pr \{ Y = j \} = \sum_{i=1}^n \Pr \{ X = i \& Y = j \}
= \sum_{i=1}^n \th_{ij} .
\ed
The vectors $\bmu \in \Sm_n , \bnu \in \Sm_m$ defined by
\bd
\mu_i = \Pr \{ X = i \} , \nu_j = \Pr \{ Y = j \} 
\ed
are called the {\bf marginal distributions} associated with $\bth$.
The ratio
\bd
\frac{ \th_{ij} }{ \nu_j } =: \Pr \{ X = i | Y = j \} 
\ed
is called the {\bf conditional probability} that $X = i$ given that $Y = j$.
Similarly,
\bd
\frac{ \th_{ij} }{ \mu_i } =: \Pr \{ Y = j | X = i \}
\ed
is the conditional probability that $Y = j$ given that $X = i$.
For each fixed $j \in \Bbb$, the vector of conditional probabilities of $X$
given that $Y = j$ belongs to $\Sm_n$, and therefore the quantity
$H(X | Y = j)$ is well-defined for each $j \in \Bbb$.
Finally, the convex combination
\be\label{eq:24d}
\sum_{j=1}^m \Pr \{ Y = j \} \cdot H(X | Y = j) =: H(X|Y)
\ee
is called the {\bf conditional entropy} of $X$ given $Y$.
It can be verified that
\be\label{eq:24e}
H(X|Y) = H(X,Y) - H(Y) = H(\bth) - H(\bnu) ,
\ee
and similarly
\bd
H(Y|X) = H(X,Y) - H(X) = H(\bth) - H(\bmu) .
\ed
respectively.
See \cite[Theorem 2.2.1]{CT06}, where (\ref{eq:24e}) is referred to
as the ``chain rule.''

\subsection{Stationary Stochastic Processes}\label{ssec:23}

Until now we have discussed either individual random variables, or
pairs of random variables.
These ideas can be extended in a straight-forward manner to
stationary stochastic processes of the form $\{ X_t \}_{t \geq 0}$ where
each $X_t$ assumes values in a finite set $\Abb$.
For each integer $t \geq 1$, we can compute the conditional entropy
$H(X_t | X_0^{t-1})$, which is the entropy of the ``current'' random
variable $X_t$ given the ``past'' $X_0^{t-1}$.
Because the process is stationary, it readily follows that
\bd
H(X_{t+\t} | X_\t^{t+\t-1} ) = H(X_t | X_0^{t-1}) ,
\fa \t \geq 0 , \fa t \geq 1.
\ed
Moreover, conditioning on more variables cannot increase entropy.
Therefore the quantity $H(X_t | X_0^{t-1})$ is nonincreasing as a
function of $t$.
Since it is clearly bounded below by zero, it follows that
there exists a constant $c$ such that
\bd
H(X_t | X_0^{t-1}) \ap c \mbox{ as } t \ap \infty .
\ed
This constant $c$ is referred to as the {\bf entropy of the process}.

Note that if the process $\{ X_t \}$ is i.i.d., then $X_t$
is independent of $X_0^{t-1}$, whence
\bd
H(X_t | X_0^{t-1}) = H(X_t) = H(X_1) , \fa t \geq 1 .
\ed
Therefore the entropy of an i.i.d.\ process is the same as the entropy
of its one-dimensional marginal.

\begin{theorem}\label{thm:21}
Suppose $\{ X_t \}$ is a Markov chain assuming values in a finite set $\Abb$,
and define $\bmu,\bmubar$ as in (\ref{eq:22a}) and (\ref{eq:22b}) 
respectively.
Then the entropy of the Markov process is given by 
\be\label{eq:24a}
c = H(\bnu) - H(\bnubar) = \sum_{i=1}^n \mubar_i H(\a^i) , 
\ee
where $\a^i$ denotes the $i$-th row of the state transition matrix $A$.
\end{theorem}

{\bf Proof:}
If $\{ X_t \}$ is a Markov process, then whenever $t \geq 2$, it follows that
\beq 
H(X_t | X_0^{t-1} ) & = & H(X_t | X_{t-1} ) = H(X_1 | X_0) \nonumber \\
& = & H(X_1,X_0) - H(X_0) , \nonumber
\eeq
where we use the stationarity of the process.
Therefore the sequence of numbers $\{ H(X_t | X_0^{t-1} ) \}$ converges
in one step to the limit $H(X_1,X_0) - H(X_0)$,
which is the entropy of a Markov process.
By definition $\bmu$ is the probability distribution of the pair $(X_0,X_1)$,
while $\bmubar$ is the probability distribution of $X_0$.
Therefore the entropy of the Markov chain is $H(\bnu) - H(\bnubar)$,
which establishes the first equality in (\ref{eq:24a}).
The second equality in (\ref{eq:24a}) follows readily from
(\ref{eq:24d}) and (\ref{eq:24e}).
We have
\beq
H(X_1 | X_0) & = & \sum_{i=1}^n \Pr \{ X_0 = i \} \cdot
H( \Pr \{ X_1 | X_0 = i \} ) \nonumber \\
& = & \sum_{i=1}^n \mubar_i H(\a^i) , \nonumber
\eeq
which is the desired conclusion.  $\halmos$


\begin{corollary}\label{corr:21}
Suppose $\{ X_t \}$ is an $s$-step Markov process, and define the
distributions $\bmu$ and $\bmubar$ as in (\ref{eq:22c}) and (\ref{eq:22d}) 
respectively.
Then the entropy of the process is given by
\be\label{eq:24b}
c = H( \bmu ) - H( \bmubar )
= \sum_{\i \in \Abb^s} \mubar_\i \sum_{j = 1}^n \frac { \mu_{\i j} }{ \mubar_\i}
\log \frac { \mu_{\i j} }{ \mubar_\i} .
\ee
\end{corollary}

The proof is very similar to that of Theorem \ref{thm:21} and is omitted.

Inspired by this discussion we now introduce a couple of definitions
that will play a direct role in the large deviation theory for Markov chains.
Suppose $\bmu \in \M_s(\Abb^k)$.
We define $\bmubar \in \M_s(\Abb^{k-1})$ by
\be\label{eq:24c}
\mubar_\i := \sum_{j \in \Abb} \mu_{\i j} , \fa \i \in \Abb^{k-1} .
\ee
The overbar serves to remind us that $\bmubar$ is ``reduced by one dimension''
from $\bmu$.
The symbol $\bnubar$ is defined similarly.

\begin{definition}\label{def:2}
Suppose $\bnu , \bmu \in \M_s(\Abb^k)$ for some integer $k \geq 2$.
Then
\be\label{eq:25}
H_c(\bmu) := H( \bmu ) - H ( \bmubar ) 
\ee
is called the {\bf conditional entropy} of $\bmu$, while
\be\label{eq:26}
D_c( \bnu \nm \bmu ) := D( \bnu \nm \bmu ) - D( \bnubar \nm \bmubar )
\ee
is called the {\bf conditional relative entropy} between $\bnu$ and $\bmu$.
\end{definition}

It is easy to show that both $H_c(\cdot)$ and $D_c(\cdot \nm \cdot)$ are 
nonnegative-valued functions.

The above definition of the entropy of a Markov process was introduced
by Shannon in \cite{Shannon48}.
He defined it as the limit of $H(X_t | X_0^{t-1})$ as $t \ap \infty$,
but did not give a closed form formula for the limit.
It was not until 2004 that an explicit expression was given in \cite{RAC04}.
It might be added that the derivation is so obvious that one wonders why
it took so long!

The remainder of this section consists of a digression to introduce
a very famous theorem known as the Shannon-McMillan-Breiman (SMB) theorem
(also known as the asymptotic equipartition property or AEP),
even though it is not directly relevant to the rest of the paper.
Let $\Abb$ be a finite set as above, and let $\Nbb$ denote the set of
natural numbers; then $\Abb^\Nbb$ denotes the set of all sequences
$\{ x_i \}_{i \geq 1}$ where each $x_i \in \Abb$.
Let $\{ X_t \}$ be a stationary ergodic process assuming values in $\Abb$.
Now define a sequence of functions (a stochastic process)
$h_l: \Abb^\Nbb \ap \R$ as follows:
For each integer $l$ and each $\x \in \Abb^\Nbb$,
\bd
h_l(\x) := - \frac{1}{l} \log \Pr \{ X_1^l = \x_1^l \} .
\ed
In other words, for a given infinite sequence $\x$, the value $h_l(\x)$
is just the negative log-likelihood that the first $l$ values of
the stochastic process $\{ X_t \}$ equal the first $l$ symbols in the
sequence $\x$.

\begin{theorem}\label{thm:SMB}
(Shannon-McMillan-Breiman Theorem)
The sequence of functions $\{ h_l(\cdot) \}$ converges in probability
and almost surely to the constant $c$ where $c$ is the entropy of the
process.
\end{theorem}

The convergence of the sequence of {\it numbers\/} $\{ H(X_t | X_0^{t-1}) \}$
as $t \ap \infty$ is a ready consequence of the fact that the sequence
is nonincreasing and is bounded below.
This limit exists for {\it every\/} stochastic process, whether ergoodic or not.
The SMB theorem states something far deeper.
The convergence in probability of the process $\{ h_l \}$ can be informally
stated as ``almost all sample paths have the same likelihood'', and is
also known as the AEP.
This theorem was first established for Markov chains by Shannon 
\cite{Shannon48}, and then extended to stationary ergodic processes by
McMillan \cite{McMillan53}, who established that the sequence of
functions $\{ h_l(\cdot) \}$ converges in probability to the constant
function $c$ where $c$ is the entropy of the process.
Finally Breiman \cite{Breiman57} showed that the sequence $\{ h_l(\cdot) \}$
also converges almost surely to the constant function $c$.

\section{Definition of the Large Deviation Property}\label{sec:3}

In this section we give a brief overview of some of the relevant
definitions from large deviation theory.
The primary references are \cite{FK06,DZ98}; specific theorem and/or
page numbers are given when appropriate.

Suppose $\{ X_t \}$ is a stationary stochastic process assuming values
in a finite set $\Abb$, and let $\Pt_\mu$ denote the law of the process;
note that $\Pt_\mu$ is a probability measure on the infinite cartesian product
$\Abb^\infty$, which can be identified with the set of all sequences taking
values in $\Abb$.
Suppose we observe a sample path $\x \in \Abb^\infty$ of this process.
Based on the first $l$ samples $x_1^l$, we can construct empirical
estimates of the statistics of the process.
For instance, we can compute an approximation $\bmuhat_k(x_1^l)$ to the
joint distribution of $k$-tuples $X_1^k$, for some fixed integer $k$.
If the process is Markov, then as discussed in Section \ref{ssec:21} we should
estimate the doublet frequency vector.
More generally, if the process is $s$-step Markov, then we should estimate
the joint distribution of $(s+1)$-tuples.
Of course, if the process is i.i.d., then we should estimate the 
one-dimensional marginal distribution of the process.

Since the computed empirical distribution $\bmuhat_k(x_1^l)$ is based on
the random sample $x_1^l$, one can think of $\{ \bmuhat_k(x_1^l) \}_{l \geq 1}$
as a stochastic process taking values in the simplex $\Sm_{n^k}$.
As an aside, until now researchers have not always ensured
that the empirically estimated distribution $\bmuhat_k(x_1^l)$ is stationary,
even though the ``true but unknown'' joint distribution of $k$-tuples
is surely stationary if $k \geq 2$.
In other words, depending on the method used, $\bmuhat_k(x_1^l)$
may belong only to $\Sm_{n^k}$ and not necessarily to $\M_s(\Abb^k)$,
even though the true distribution $\bmu$ belongs to $\M_s(\Abb^k)$.
In the present paper, we take special care to ensure that all empirical
estimates are stationary.

Since $\M_s(\Abb^k)$ is a subset of $\Sm_{n^k}$ which is in turn a subset
of $\R^{n^k}$, we can equip both $\Sm_{n^k}$ and $\R^{n^k}$ with the
$\ell_1$-norm, usually referred to as the {\bf total variation
metric}.\footnote{Actually, the total variation metric between two
probability distributions $\bnu$ and $\bmu$ equals {\em one half\/}
of the $\ell_1$-norm $\nm \bnu - \bmu \nm_1$.}
Then phrases such as open sets and closed sets in $\M_s(\Abb^k)$ and
$\Sm_{n^k}$ have an unambiguous meaning.

Next we define the large deviation property and the large deviation
rate function for the $\Sm_{n^k}$-valued stochastic process
$\{ \bmuhat_k(x_1^l) \}$.

\begin{definition}\label{def:LDP}
The process $\{ \bmuhat_k(x_1^l) \}$ is said to
satisfy a {\bf large deviation property (LDP) with rate function
$I : \Sm_{n^k} \ap \R \cup \{ \infty \}$}
if the following two statements hold:
\ben
\item
For each open set $\GA \seq \Sm_{n^k}$, we have that
\be\label{eq:31}
- \inf_{\bnu \in \GA} I(\bnu) \leq
\liminf_{\lai} \frac{1}{l} \log \Pt_\mu \{ \bmuhat_k(x_1^l) \in \GA \} .
\ee
\item For each closed set $\GA \seq \Sm_{n^k}$, we have that
\be\label{eq:32}
\limsup_{\lai} \frac{1}{l} \log \Pt_\mu \{ \bmuhat_k(x_1^l) \in \GA \} \leq
- \inf_{\bnu \in \GA} I(\bnu) .
\ee
\een
\end{definition}

As shown in \cite[Section 3.1]{FK06}, one can assume without loss
of generality that the function $I$ is lower semicontinuous by replacing
$I$ by its lower semicontinuous relaxation if necessary.
Moreover, there is at most one lower semicontinuous function $I$
that satisfies (\ref{eq:31}) and (\ref{eq:32}).
Thus, if at all a LDP holds, then the corresponding lower semi-continuous
rate function is uniquely determined.
Moreover, if we assume (without loss of generality) that the rate
function $I$ is lower semicontinuous, then the two equations
(\ref{eq:31}) and (\ref{eq:32}) can be combined into the following
single statement:
Let $\GA \seq \Sm_{n^k}$ be any Borel set of probability distributions.
Then
\beq
- \inf_{\bnu \in \GA^o} I(\bnu)
& \leq & \liminf_{\lai} \frac{1}{l} \log \Pt_\mu \{ \bmuhat_k(x_1^l) \in \GA \}
\nonumber \\
& \leq & \limsup_{\lai} \frac{1}{l} \log \Pt_\mu \{ \bmuhat_k(x_1^l) \in \GA \}
\nonumber \\
& \leq & - \inf_{\bnu \in \bar{\GA} } I(\bnu) . \label{eq:33}
\eeq
In the above formula, $\GA^o$ denotes the interior of $\GA$,
while $\bar{\GA}$ denotes the closure of $\GA$, both in the total
variation metric.
In particular, suppose $I$ is continuous (not merely lower semicontinuous),
and that $\GA$ does not have any isolated points; that is,
$\GA \seq \overline{ \GA^o }$.
Then the two extreme infima in (\ref{eq:33})
are equal, which means that the liminf, limsup,
and limit all coincide, and we can conclude that
\bd
\lim_{\lai} \frac{1}{l} \log \Pt_\mu \{ \bmuhat_k(x_1^l) \in \GA \}
= - \inf_{\bnu \in \bar{\GA} } I(\bnu) .
\ed
This precise estimate of the rate of convergence of empirical estimates
to their true values is what gives large deviation theory its power
and appeal.

\section{Summary of Known Rate Functions}\label{sec:4}

In this section we summarize several known results on rate functions
for i.i.d.\ processes and for Markov chains.
Recall that the main objective of the paper is to present elementary
proofs of most of these results.

\subsection{Sanov's Theorem for Finite Alphabets}\label{ssec:41}

In this subsection we state a well-known result known as Sanov's theorem,
which gives the rate function for i.i.d.\ processes.

Suppose $\{ X_t \}$ is an i.i.d.\ process assuming values in a finite
set $\Abb$ with the one-dimensional marginal distribution $\bmu \in \Sm_n$.
Based on an observation $x_1^l$, we construct the empirical measure
$\bphi(x_1^l) \in \Sm_n$ as
\be\label{eq:34}
\phi_j(x_1^l) := \frac{1}{l} \sum_{t=1}^l I_{ \{ x_t = j \} } ,
\fa j \in \Abb ,
\ee
where $I$ denotes the indicator function.
Therefore $\phi_j(x_1^l)$ is just the fraction of the first $l$ samples
that are equal to the symbol $j$.
Note that $\bphi(x_1^l)$ is an empirical estimate of $\bmu$,
and could therefore be denoted by $\bmuhat_1(x_1^l)$.
However, to reduce notational clutter we use $\bphi(x_1^l)$ instead.
The well-known Sanov's theorem when specialized to a finite alphabet
gives an explicit formula for the rate function of the process 
$\{ \bphi(x_1^l) \}$; see for example \cite[Section 2.1.1]{DZ98}.
 
\begin{theorem}\label{thm:31}
The $\Sm_n$-valued process
$\{ \bphi(x_1^l) \}$ satisfies the LDP with the rate function
$I(\bnu) = D( \bnu \nm \bmu )$.
\end{theorem}

\subsection{Rate Function for Singleton Frequencies of a Markov
Chain}\label{ssec:42}

Throughout the remainder of the section, the object of study is a Markov process
$\{ X_t \}$ assuming values in a finite set $\Abb$, with doublet
frequency vector $\bmu \in \McA$, stationary distribution $\bmubar \in \Sm_n$
and state transition matrix $A \in [0,1]^{n \times n}$ defined as before by
\bd
\mubar_i = \sum_{j=1}^n \mu_{ij} ,
a_{ij} = \frac{ \mu_{ij} }{ \mubar_i } .
\ed

Next we discuss the rate function for the singleton frequency distribution
of a Markov chain, or equivalently, the stationary distribution of a
Markov chain.

Suppose $\{ X_t \}$ is a Markov process assuming values in a finite set
$\Abb$.
Given an observation $x_1^l = x_1 \ldots x_l$, we can as before form
an empirical distribution $\bphi(x_1^l) \in \Sm_n$ in analogy with
(\ref{eq:34}); that is,
\be\label{eq:34a}
\phi_j(x_1^l) := \frac{1}{l} \sum_{t=1}^l I_{ \{ x_t = j \} } ,
\fa j \in \Abb .
\ee
Thus
$\bphi$ is an approximation to the stationary distribution $\bpi$ of
the Markov chain.

The rate function for the $\Sm_n$-valued process $\{ \bphi(x_1^l) \}$
was first derived by Donsker and Varadhan in a series of papers
\cite{DV75b,DV75c,DV76}, using a characterization of the spectral radius
of a positive matrix \cite{DV75a}.
In \cite{DZ98}, an alternate derivation of these results
is given based on the G\"{a}rtner-Ellis theorem \cite{Gartner77,Ellis84}.

\begin{theorem}\label{thm:32}
Suppose $\{ X_t \}$ is a Markov process in the finite alphabet $\Abb$
with doublet frequency distribution $\bmu \in \McA$, stationary
distribution $\bmubar \in \Sm_n$, and 
state transition matrix $A \in [0,1]^{n \times n}$.
Suppose the state transition matrix $A$ of the Markov process $\{ X_t \}$
is irreducible.
Then the $\Sm_n$-valued process $\{ \bphi(x_1^l) \}$ satisfies the LDP
with the rate function
\be\label{eq:35}
I(\bphi) = \sup_{\u > \bz} \sum_{i=1}^n \phi_i \log \frac{u_i}{(\u A)_i} ,
\ee
\be\label{eq:35a}
I(\bphi) = \sup_{\u > \bz} \sum_{i=1}^n \phi_i \log \frac{u_i}{( A \u )_i} .
\ee
\end{theorem}

Note that (\ref{eq:35}) is given in Theorem 3.1.6 and (\ref{eq:35a}) is
given in Exercise 3.1.11 of \cite{DZ98}.
Here the notation $\u > \bz$ means that $u_i > 0 \fa i$.

\subsection{Rate Functions for Doublet Frequencies of a Markov
Chain}\label{ssec:43}

Given a sample path $x_1^l$ of the Markov chain,
there is more than one way to construct an empirical estimate of
the doublet frequency distribution.
The precise estimate chosen has an impact on the difficulty of the analysis.

Given the sample path $x_1^l = x_1 \ldots x_l$, one possibility is to define
\be\label{eq:36}
\th_{ij}(x_1^l) := \frac{1}{l-1} \sum_{t=1}^{l-1} I_{ \{ X_t X_{t+1} = ij \} } .
\ee
This procedure produces a vector $\bth \in \Sm_{n^2}$ which can be
interpreted as a measure on $\Abb^2$.
However, the distribution $\bth$ is {\em not stationary\/} in general.
If we define $\bthbar \in \Sm_n$ by
\bd
\thbar_i := \frac{1}{l-1} \sum_{t=1}^{l-1} I_{ \{ X_t = i \} } ,
\ed
then it is certainly true that
\bd
\thbar_i = \sum_{j=1}^n \th_{ij} .
\ed
However, in general
\bd
\sum_{j=1}^n \th_{ji} \neq \thbar_i .
\ed
Hence $\bth \in \Sm_{n^2}$ is {\em not\/} a stationary distribution
in general.
Moreover, there is no simple relationship between $\bar{\bth} \in \Sm_n$
and $\bphi \in \Sm_n$ defined in (\ref{eq:34a}).
If $x_l = x_1$ so that the sample path is a cycle, then $\bth \in \McA$,
but not in general.
Moreover, even if $x_l = x_1$ so that the sample path is a cycle,
in general $\bthbar \neq \bphi$.

The rate function for the process $\{ \bth(x_1^l) \}$
is given in \cite[Theorem 3.1.13]{DZ98}.

\begin{theorem}\label{thm:33}
Suppose $\{ X_t \}$ is a Markov process in the finite alphabet $\Abb$
with doublet frequency distribution $\bmu \in \McA$, stationary
distribution $\bmubar \in \Sm_n$, and 
state transition matrix $A \in [0,1]^{n \times n}$.
Suppose $a_{ij} > 0 \fa i, j \in \Abb$.
Then the $\Sm_{n^2}$-valued process $\{ \bth(x_1^l) \}$ satisfies the LDP
with the rate function
\be\label{eq:37}
I(\bth) =
\sum_{i \in \Abb} \bar{\th}_i \sum_{j \in \Abb} b_{ij} \log (b_{ij}/a_{ij}) ,
\fa \bth \in \McA , \\
\ee
\be\label{eq:37a}
I(\bth) = + \infty ,
 \fa \bth \not \in \McA .
\ee
Here, as before, $\bar{\bth} \in \Sm_n$ is the one-dimensional marginal
of $\bth$ and 
\bd
b_{ij} = \frac{\th_{ij}}{\bar{\th}_i} 
\ed
is the state transition matrix associated with $\bth$ 
if $\bth \in \McA$.
\end{theorem}

The proof of the above theorem as given in \cite[Theorem 3.1.13]{DZ98}
is based on the observation that
if $\{ X_t \}$ is a Markov chain, then so is the stochastic process consisting
of doublets $\{ (X_t, X_{t+1}) \}$.
By applying Theorem \ref{thm:32} to the latter Markov chain, it is possible
to compute the rate function of doublet frequencies of the original
Markov chain.
A perusal of the relevant pages shows that the proof of this theorem
using this approach is anything but simple.

If we amend slightly the manner in which
the empirical estimate is constructed so that the estimated vector is
{\em always stationary}, then the derivation of the rate function
is greatly simplified.
This is the approach adopted in \cite{Natarajan85}.

Given a sample path $x_1^l = x_1 \ldots x_l$,
we construct the empirical estimate $\bnu = \bnu(x_1^l)$ as follows:
\be\label{eq:41}
\nu_{ij}(x_1^l) := \frac{1}{l} \sum_{t=1}^l I_{ \{ x_t x_{t+1} = ij \} } ,
\fa i,j \in \Abb ,
\ee
where $x_{l+1}$ is taken as $x_1$.
If we compare (\ref{eq:41}) with (\ref{eq:35}), we see that we have
in effect augmented the original sample path $x_1^l$ by adding
a ``ghost'' transition
from $x_l$ back to $x_1$ so as to create a cycle, and used this artificial 
sample path of length $l+1$ to construct the empirical estimate.
The advantage of doing so is that the resulting vector $\bnu$ is 
{\em always stationary}, unlike $\bth$ in (\ref{eq:35}) which may not be
stationary in general.
Moreover, the one-dimensional marginal of $\bnu$ is precisely
$\bphi$ defined in (\ref{eq:34}).
This is established next.

Recall the earlier definition, reproduced here for convenience:
\bd
\phi_j(x_1^l) := \frac{1}{l} \sum_{t=1}^l I_{ \{ x_t = j \} } ,
\fa j \in \Abb .
\ed
Now it is claimed that if $\nu_{ij}(x_1^l)$ is defined as in (\ref{eq:41}), 
then
\be\label{eq:43}
\sum_{i \in \Abb} \nu_{ji} = \phi_j , \fa j \in \Abb ,
\sum_{j \in \Abb} \nu_{ji} = \phi_i , \fa i \in \Abb ,
\ee
thus showing that $\bnu$ is stationary.
To establish (\ref{eq:43}),
observe that the quantity $\sum_{i \in \Abb} \nu_{ji}$ is obtained by
counting the
number of times that the symbol $j$ occurs as the first symbol in the
sequence $x_1 x_2 , x_2 x_3, \ldots , x_{l-1} x_l,x_l x_1$,
and then dividing by $l$.
Similarly the quantity $\sum_{j \in \Abb} \nu_{ji}$ is obtained by counting
the number of times that $i$ occurs as the second symbol 
in $x_1 x_2 , x_2 x_3, \ldots , x_{l-1} x_l,x_l x_1$, and then dividing by $l$.
Now, for $2 \leq t \leq l-1$, $x_t$ is the first symbol in $x_t x_{t+1}$,
and the second symbol in $x_{t-1} x_t$.
Next, $x_1$ is the first symbol in $x_1 x_2$ and the second symbol in the
ghost transition $x_l x_1$.
Similarly $x_l$ is the second symbol in $x_{l-1} x_l$ and the first symbol
in the ghost transition $x_l x_1$.
This establishes (\ref{eq:43}) and shows that $\bnu$ is stationary.
Moreover, its one-dimensional marginal is $\bphi$ as defined in (\ref{eq:34}).

The rate function for the alternative, and stationary, estimate
$\{ \bnu(x_1^l) \}$ is given next; see \cite[Theorem 1]{Natarajan85}.

\begin{theorem}\label{thm:34}
Suppose $\{ X_t \}$ is a Markov process in the finite alphabet $\Abb$
with doublet frequency distribution $\bmu \in \McA$, stationary
distribution $\bmubar \in \Sm_n$, and 
state transition matrix $A \in [0,1]^{n \times n}$.
Suppose $\mu_{ij} > 0 \fa i, j \in \Abb$.
Then the $\McA$-valued process $\{ \bnu(x_1^l) \}$ satisfies the LDP
with the rate function
\be\label{eq:43a}
I(\bnu) := 
\sum_{i \in \Abb} \bar{\nu}_i \sum_{j \in \Abb} c_{ij} \log (c_{ij}/a_{ij})
\ee
where $\bar{\bnu} \in \Sm_n$ is the one-dimensional marginal
of $\bnu$ defined in analogy with (\ref{eq:24c}), and
\bd
c_{ij} = \frac{\nu_{ij}}{\bar{\nu}_i} .
\ed
\end{theorem}

{\bf Remarks:}
\ben
\item It is clear that, for a given finite sample path $x_1^l$, the two
empirical doublet distributions $\bth(x_1^l)$ and $\bnu(x_1^l)$ are
in general different.
For one thing, $\bnu(x_1^l)$ is always stationary, whereas $\bth(x_1^l)$
need not be stationary in general.
Now suppose $x_l = x_1$ so that the sample path $x_1^l$ is a cycle,
and as a result $\bth(x_1^l)$ is also stationary.
Even in this case, the two empirical doublet distributions $\bth(x_1^l)$
and $\bnu(x_1^l)$ are in general different.
This is because $\bth(x_1^l)$ is based on the original sample path of
length $l$, whereas $\bnu(x_1^l)$ is based on an augmented
sample path of length $l+1$.
\item
One consequence of this is that, even when $\bth(x_1^l)$ is stationary,
in general its one-dimensional marginal bears no relationship to the empirical
estimate $\bphi(x_1^l)$ defined in (\ref{eq:34}), which is a natural
definition for both i.i.d.\ processes as well as for Markov processes.
In contrast, not only is $\bnu(x_1^l)$ always guaranteed to be stationary,
but its one-dimensional marginal is indeed $\bphi(x_1^l)$.
\item
It is equally clear that the two rate functions in (\ref{eq:37}) and
(\ref{eq:43a}) are identical.
Thus introducing the `ghost' transition from $x_l$ to $x_1$ to ensure
that the estimated doublet frequency vector is stationary does not
affect the rate function.
However, as we shall see below, this makes the {\em derivation of the rate
function\/} very easy.
\item In (\ref{eq:37}) we are obliged to specify the rate function over
{\em all\/} of $\Sm_{n^2}$ because in general the estimate $\bth(x_1^l)$
need not be stationary.
In contrast, the estimate $\bnu(x_1^l)$ is always
guaranteed to be stationary.
Hence it is necessary to define the rate function only over $\McA$.
\item The hypothesis in Theorem \ref{thm:34}
that $\bmu$ is strictly positive is clearly
equivalent to the hypothesis of Theorem \ref{thm:33}, namely that
$a_{ij} > 0 \fa i, j$.
This is because $a_{ij} = \mu_{ij}/\mubar_i$.
\een

The construction in (\ref{eq:41}) is intended to generate a stationary
estimate belonging to $\McA$ for (one-step) Markov chains.
By following a similar philosophy, it is possible to generate a stationary
estimate belonging to $\McAs$ for $s$-step Markov chains
Suppose $\{ X_t \}$ is an $s$-step Markov process, so that
\bd
E \{ X_t | X^{t-1}_0 \} = E \{ X_t | X^{t-1}_{t-s} \}
\fa t .
\ed
Then the process is completely characterized by the $(s+1)$-tuple
frequency vector
\bd
\mu_\i := \Pr \{ X_t^{t+s} = \i \} , \fa \i \in \Abb^{s+1} .
\ed
Note that the frequency vector $\bmu$ is stationary and thus belongs to $\McAs$.
Since an $s$-step Markov process over $\Abb$ can be viewed as a conventional
(one-step) Markov process over the state space $\Abb^s$, we can identify
the stationary distribution
\bd
\mubar_\i := \sum_{j \in \Abb} \mu_{i \j} =
\sum_{j \in \Abb} \mu_{j \i} , \fa \i \in \Abb^s ,
\ed
while the transition probabilities are given by
\bd
\Pr \{ X_t = j | X_{t-s}^{t-1} = \i \} = \frac{ \mu_{\i j}}{\mubar_\i} .
\ed

Suppose $x_1^l$ is a sample path of length $l$ of an $s$-step Markov chain.
To construct a {\em stationary\/} empirical measure on the basis of
this sample path, we define the augmented sample path $\tilde{x}_1^l
:= x_1 \ldots x_l x_1 \ldots x_s = x_1^l \cdot x_1^s \in \Abb^{l+s}$.
Here the symbol $\cdot$ denotes the concatenation of two strings.
The above augmentation is the $s$-step generalization of adding a single
ghost transition from $x_l$ to $x_1$ in the case of one-step Markov chains.
In this case we are adding $s$ ghost transitions.
Then we define
\be\label{eq:61}
\nu_\i := \frac{1}{l} \sum_{t=1}^l I_{ \{ x_t^{t+s} = \i \} } ,
\fa \i \in \Abb^{s+1} .
\ee
Note that (\ref{eq:61}) reduces to (\ref{eq:41}) if $s = 1$.
Then the resulting empirical measure $\bnu(x_1^l)$ belongs to $\McAs$.
For this empirical measure, there aren't too many results that are
explicitly stated concerning its LDP rate function.
However, as we shall see below, the method of proof used here extends in
a straight-forward manner to multi-step Markov chains, with the consequence
that obtaining LDP rate functions for multi-step Markov chains is not
any more difficult than doing so for single-step Markov chains.

\section{The Method of Types}\label{sec:4a}

For the purposes of the present paper, the method of types can be described
as follows:
For a given integer $l$, the sample space of all possible sample paths is
clearly $\Abb^l$.
With each sample path $x_1^l \in \Abb^l$, we associate a corresponding
empirical distribution $\bnu(x_1^l) \in \McA$ as in (\ref{eq:41})
in the case of one-step Markov chains, and $\bnu(x_1^l) \in \McAs$
as in (\ref{eq:61}) in the case of $s$-step Markov chains.
The essence of the method of types is to consider two sample paths
$x_1^l, y_1^l \in \Abb^l$ to be equivalent if they lead to the
same empirical distribution, that is, if $\bnu(x_1^l) = \bnu(y_1^l)$.
(It is easy to see that this does indeed define an equivalence relation
on $\Abb^l$.)
Let $\E(l,n,2)$ denote the set of all possible empirical distributions
$\bnu(x_1^l) \in \McA$ that can be generated using (\ref{eq:41})
as $x_1^l$ varies over $\Abb^l$.
Note that the set $\E(l,n,2)$ depends on both $l$, the length of the sample
path, as well as $n$, the cardinality of the underlying state space $\Abb$.
The symbol $2$ serves to remind us that we are computing doublet frequencies.
More generally, let $\E(l,n,s+1)$ denote the set of all possible
empirical distributions $\bnu(x_1^l) \in \McAs$
that can be generated using (\ref{eq:61})
as $x_1^l$ varies over $\Abb^l$.
For each $\bzeta \in \E(l,n,2)$, define
\bd
T(\bzeta,l,2) := \{ x_1^l \in \Abb^l : \bnu(x_1^l) = \bzeta \} .
\ed
Thus $T(\bzeta,l,2)$ is the set of all sample paths of length $l$ that
generate the empirical distribution $\bzeta$.
Then $T(\bzeta,l,2) \seq \Abb^l$ is called the {\bf type class} of $\bzeta$.
Similarly, for each $\bzeta \in \E(l,n,s+1)$, define
\bd
T(\bzeta,l,s+1) := \{ x_1^l \in \Abb^l : \bnu(x_1^l) = \bzeta \} .
\ed
Strictly speaking, we should use different symbols in (\ref{eq:41}) and in
(\ref{eq:61}).
But the context should make it clear which definition is being used,
and in any case (\ref{eq:61}) reduces to (\ref{eq:41}) when $s = 1$.

The derivation of the rate function for $\bnu(x_1^l)$,
based on the method of types,
consists in addressing the following questions:
\bit
\item What is the cardinality of $\E(l,n,s+1)$?
In other words, how many distinct empirical measures $\bnu(x_1^l) \in \McAs$
can be generated as $x_1^l$ varies over $\Abb^l$?
\item For a given $\bzeta \in \E(l,n,s+1)$,
what is the cardinality of the associated type class $T(\bzeta,l,s+1)$?
In other words, how many different sample paths of length $l$ can
produce a given $\bzeta \in \E(l,n,s+1)$?
Here we require both upper as well as lower bounds on $| T(\bzeta,l,s+1) |$.
\item What is the (log) likelihood of each sample path in $T(\bzeta,l,s+1)$,
and how is it related to $\bzeta$?
\eit


Note that each $\bnu(x_1^l) \in \E(l,n,2)$
is of the form $\nu_{ij} = l_{ij}/l$ for some integer $l_{ij}$.
Moreover, the corresponding reduced distribution $\bnubar$ over $\Abb$
is given by $\nubar_i = \lbar_i/l$ where
\bd
\lbar_i = \sum_{j=1}^n l_{ij} = \sum_{j = 1}^n l_{ji} , \fa i .
\ed
Throughout the proof, $l$ denotes the length of the sample path and
$l_{ij}, \lbar_i$ denote the integers defined above.
More generally,
suppose $\bnu \in \E(l,n,s+1)$, and for each string $\j \in \Abb^{s+1}$,
let $l_{\j}$ denote the integer $l \nu_{\j}$, and let $\lbar_\i$
denote the integer $l \nubar_\i$.

\begin{lemma}\label{lemma:41}
For each $l$, we have
\be\label{eq:47a}
| \E(l,n,s+1) | \leq (l+1)^{s+1} . 
\ee
\end{lemma}

{\bf Proof:}
The proof closely follows that in \cite[Section 2.1.1]{DZ98}.
Suppose $\bzeta \in \E(l,n,s+1)$.
Then each component $\zeta_\j$ has $l+1$ possible values, namely
$0/l , 1/l, \ldots , (l-1)/l , l/l$, and there are $s+1$ such
components.
Thus the maximum number of possible vectors in $\E(l,n,s+1)$ is $(l+1)^{s+1}$.
$\halmos$ 

\begin{lemma}\label{lemma:41a}
For all $\bpsi \in \McAs$, we have that
\be\label{eq:47b}
\min_{ \bnu \in \E(l,n,s+1) } \nm \bpsi - \bnu \nm_1 \leq 
\frac{ 2 (s+3) n^{s+1} }{l} .
\ee
\end{lemma}

The proof of this lemma is beyond the scope of the paper,
and the reader is referred to \cite[Lemma 5.2]{CG05}.

\begin{lemma}\label{lemma:42}
Suppose $\bzeta \in \E(l,n,s+1)$.
Then the cardinality of the type class $T(\bzeta,l,2)$ is bounded by
\be\label{eq:47c}
(2l)^{-n^{s+1}} e^{lD(\bzeta)} \leq | T(\bzeta,l,s+1)| \leq l e^{lD(\bzeta)} 
\fa l \geq n .
\ee
\end{lemma}

Because of the complexity of the proof, Lemma \ref{lemma:42} is
stated as Theorem \ref{thm:A1} and is proved in Appendix A.


\section{Rate Function for Doublet Frequency Vector of a Markov Chain:
Statement and Proof}\label{sec:5}

The next two sections form the heart of the paper.
In this section we present the rate function for the process
$\{ \bnu(x_1^l) \}$ 
as the conditional relative entropy between the ``empirical'' and ``true''
doublet frequencies, and then give a proof based on the method of types.
This approach was introduced by Csisz\'{a}r and his coworkers; see the book by
Csisz\'{a}r and K\"{o}rner \cite{CK81} for a detailed exposition, and
the survey paper by Csisz\'{a}r \cite{Csiszar98} for a ``condensed'' version.
The original method of types was for i.i.d.\ processes.
In \cite{DLS81}, some bounds on the cardinality of type classes were
derived for Markov processes, and in \cite{Natarajan85} these bounds
were used to derive the LDP rate function.
In \cite{CCC87}, a similar analysis is performed from the standpoint of
source coding.
In the recent book \cite{Hollander00}, many of these bounds are reproduced.
Indeed, the main arguments in Chapter IV of \cite{Hollander00} are virtually
identical to the treatment here.
However, the proof of the LDP given in \cite{Hollander00}, Section IV.2,
makes use of Varadhan's lemma, a very powerful and general result, which
is far more powerful than the situation demands.
In contrast, in the proof below we pull together relevant estimates
from various sources and avoid making use of any advanced concepts.

By stating the LDP rate function as the conditional relative entropy,
we are able to extend the rate function in a transparent manner to
multi-step Markov chains.
Moreover, with some technical conditions, the LDP rate function can also
be extended to stationary processes with ``fading memory.''
Note that a multi-step Markov chain can be thought of as a stationary
process with finite memory.
Indeed, in the paper \cite{CG05}, this is the approach adopted.
However, in the present paper the focus is restricted only to (finite memory)
Markov processes.
In the proof of Theorem \ref{thm:41}
we consolidate in one place various arguments and bounds scattered
over the literature, and make them available in one place.
As a result, the proof given here is both complete and elementary,
avoiding all advanced concepts.

\begin{theorem}\label{thm:41}
Suppose $\{ X_t \}$ is a Markov process in the finite alphabet $\Abb$
with doublet frequency distribution $\bmu \in \McA$.
Suppose $\mu_{ij} > 0 \fa i, j \in \Abb$.
Then the $\McA$-valued process $\{ \bnu(x_1^l) \}$ satisfies the LDP
with the rate function
\be\label{eq:47}
I(\bnu) := D_c( \bnu \nm \bmu )
\ee
\end{theorem}

Before giving a proof of Theorem \ref{thm:41},
we show that the two quantities in
(\ref{eq:47}) and (\ref{eq:43a}) are the same.
Thus Theorems \ref{thm:34} and \ref{thm:41} are the same.
To establish the equivalence of the two formulas in (\ref{eq:47}) and
(\ref{eq:43a}), let us recall the notation in Theorem \ref{thm:34},
namely $a_{ij} = \mu_{ij}/\mubar_i$, and $c_{ij} = \nu_{ij} / \nubar_i$.
Therefore
\bd
D( \bnu \nm \bmu ) = \sum_{i \in \Abb} \sum_{j \in \Abb}
\nu_{ij} \log \frac{ \nu_{ij} }{ \mu_{ij} } ,
\ed
while
\beq
D( \bnubar \nm \bmubar ) & = & 
\sum_{i \in \Abb} \nubar_i \log \frac { \nubar_i }{ \mubar_i } \nonumber \\
& = & \sum_{i \in \Abb} \left[ \sum_{j \in \Abb} \nu_{ij} \right]
\log \frac { \nubar_i }{ \mubar_i } . \nonumber 
\eeq
Therefore
\beq
D_c( \bnu \nm \bmu ) & = & D( \bnu \nm \bmu ) - D( \bnubar \nm \bmubar )
\nonumber \\
& = &  \sum_{i \in \Abb} \sum_{j \in \Abb}
\nu_{ij} \log \frac{ \nu_{ij} / \nubar_i }{ \mu_{ij} / \mubar_i } \nonumber \\
& = & \sum_{i \in \Abb} \nubar_i \sum_{j \in \Abb}
c_{ij} \log \frac{ c_{ij} }{ a_{ij} } , \label{eq:47d}
\eeq
which is the same as (\ref{eq:43a}).

{\bf Proof of Theorem \ref{thm:41}:}
Suppose we have a sample path $x_1^l$.
Let us compute its likelihood in terms of the properties of the corresponding
empirical distribution $\bnu(x_1^l)$.
We have
\beq
\Pr \{ X_1^l = x_1^l \} & = & \Pr \{ X_1 = x_1 \} \nonumber \\
& \cdot & \prod_{t=1}^{l-1} \Pr \{ X_{t+1} = x_{t+1} | X_t = x_t \} .
\nonumber
\eeq
Hence\footnote{In the interests of clarity, in the proof we write
$\mu(x_s x_t)$ instead of $\mu_{x_s x_t}$,
and $\mubar(x_t)$ instead of $\mubar_{x_t}$.
However, we continue to use the subscript notation if the arguments are
simple indices such as $i$ and $j$.}
\beq
\Pr \{ X_1^l = x_1^l \} & = & \mubar(x_1) \cdot
\prod_{t=1}^{l-1} \frac{ \mu( x_t x_{t+1} ) }{ \mubar(x_t) } \nonumber \\
& = & \mubar(x_1) \cdot
\prod_{t=1}^l \frac{ \mu( x_t x_{t+1} ) }{ \mubar(x_t) } 
\cdot \frac{ \mubar (x_l) }{ \mu(x_l x_1) } \nonumber \\
& = & \frac{ \mubar (x_1) \mubar (x_l) }{ \mu(x_l x_1) } \cdot
\prod_{t=1}^l \frac{ \mu( x_t x_{t+1} ) }{ \mubar(x_t) } , \nonumber
\eeq
where as before we take $x_{l+1} = x_1$.
Now, since $\mu_{ij} > 0$ for all $i,j$, there exist
constants $\underline{a}, \bar{a}, \underline{b}, \bar{b}$ such that
\bd
0 < \underline{a} \leq \mubar_i \leq \bar{a} \fa i \in \Abb ,
0 < \underline{b} \leq \mu_{ij} \leq \bar{b} \fa i, j \in \Abb .
\ed
Now define constants $\underline{c}$ and $\bar{c}$ by
\bd
\underline{c} = \log \frac{ \underline{a}^2 }{ \bar{b} }
= 2 \log \underline{a} - \log \bar{b} ,
\ed
\bd
\bar{c} = \log \frac{ \bar{a}^2 }{ \underline{b} }
= 2 \log \bar{a} - \log \underline{b} .
\ed
Then it is obvious that
\bd
\underline{c} \leq \log \frac{ \mubar_i \mubar_j }{ \mu_{ij} } 
\leq \bar{c}, \fa i, j.
\ed
Of course the constants $\underline{c}$ and $\bar{c}$
depend on $\bmu$, but the point is that they
do not depend on the empirical measure $\bnu(x_1^l)$.
Therefore it follows that
\beq
\log \prod_{t=1}^l \frac{ \mu( x_t x_{t+1} ) }{ \mubar(x_t) } + \underline{c}
& \leq & \log \Pr \{ X_1^l = x_1^l \} \nonumber \\
& \leq & \log \prod_{t=1}^l \frac{ \mu( x_t x_{t+1} ) }{ \mubar(x_t) } + \bar{c} . 
\nonumber
\eeq

Next we examine the logarithm of the product term above.
We have
\bd
\log \left[ \prod_{t=1}^l \frac{ \mu( x_t x_{t+1} ) }{ \mubar(x_t) } \right]
= \sum_{t=1}^l [ \log \mu( x_t x_{t+1} ) - \log \mubar(x_t) ] .
\ed
When we do the above summation, we observe that the event $x_t x_{t+1} = ij$
occurs exactly $l_{ij} = l [\bnu (x_1^l)]_{ij}$ times, while the event
$x_t = i$ occurs exactly $\bar{l}_i = l [\bnubar(x_1^l)]_i$ times.
Therefore, instead of summing over $t$, we can sum over all $i$ and $j$.
This gives
\beq
\log \left[ \prod_{t=1}^l \frac{ \mu( x_t x_{t+1} ) }{ \mubar(x_t) } \right]
& = & l \sum_{i \in \Abb} \sum_{j \in \Abb} \nu_{ij} \log \mu_{ij} \nonumber \\
& - & l \sum_{i \in \Abb} \nubar_i \log \mubar_i \nonumber \\
& = & -l [ J( \bnu,\bmu ) - J( \bnubar , \bmubar )] , \nonumber
\eeq
where we write $\bnu$ and $\bnubar$ for the more precise $\bnu(x_1^l)$
and $\bnubar(x_1^l)$.
Substituting this into the previous bound leads to the following relationships:
\be\label{eq:49a}
\log \Pr \{ X_1^l = x_1^l \} \geq - l [ J( \bnu,\bmu )
- J( \bnubar , \bmubar )] + \underline{c} ,
\ee
\be\label{eq:49b}
\log \Pr \{ X_1^l = x_1^l \} \leq - l [ J( \bnu,\bmu ) 
- J( \bnubar , \bmubar )] + \bar{c} ,
\ee
where $\bnu$ is a shorthand for $\bnu(x_1^l)$.

In large deviation theory, the quantity of interest is the log of the 
likelihood that a particular empirical estimate will occur, normalized
by the length of the observation.
Accordingly, let us denote the empirical distribution generated
by a sample path as $\bnu(x_1^l)$, and for each $\bzeta \in \E(l,n,2)$ define
\bd
\d(l,\bzeta) := \frac{1}{l} \log \Pr \{ \bnu(x_1^l) = \bzeta \} .
\ed
What interests us is the behavior of $\d(l,\bzeta)$ as $\lai$.
Now we know from (\ref{eq:49a}) and (\ref{eq:49b}) 
that the normalized log likelihood of each sample path
within the type class $T(\bzeta,l,2)$ looks like
$J( \bzeta,\bmu ) - J( \bzetabar , \bmubar ) + o(1/l)$, and we know from
(\ref{eq:47c}) that $(1/l) \log |T(\bzeta,l,2)|$ looks like
$H_c(\bzeta) + o(1/l)$.
Combining these two facts leads to
\beq
\d(l,\bzeta)
& \leq & H_c(\bzeta) - J(\bzeta,\bmu) + J(\bzetabar,\bmubar) + o(1/l) \nonumber \\
& = & H(\bzeta) - H(\bzetabar) - J(\bzeta,\bmu) + J(\bzetabar,\bmubar)
+ o(1/l) \nonumber \\
& = & - D(\bzeta \nm \bmu) + D(\bzetabar \nm \bmubar)+ o(1/l) \nonumber \\
& = & - D_c(\bzeta \nm \bmu) + o(1/l) . \label{eq:410}
\eeq
Similarly we get
\be\label{eq:411}
\d(l,\bzeta)
\geq - D_c(\bzeta \nm \bmu) + o(1/l) .
\ee

The completion of the proof is based on adapting absolutely standard techniques
to the present situation, and the reader
is referred to \cite[pp.\ 16-17]{DZ98} for comparison purposes.
Let $\GA \seq \McA$ be any set of stationary distributions on $\Abb^2$.
Then
\beq
\Pr \{ \bnu(x_1^l) \in \GA \} & = & 
\sum_{ \bzeta \in \E(l,n,2) \cap \GA } \Pr \{ \bnu(x_1^l) = \bzeta \} \nonumber \\
& \leq & | \E(l,n,2)  | \sup_{ \bzeta \in \GA }
\Pr \{ \bnu(x_1^l) = \bzeta \} . \nonumber 
\eeq
Hence
\bd
\frac{1}{l} \log \Pr \{ \bnu(x_1^l) \in \GA \}
\leq \frac{1}{l} \log | \E(l,n,2) | + \sup_{\bzeta \in \GA } \d(l,\bzeta) .
\ed
Since $|\E(l,n,2)|$ is polynomial in $l$,
it follows that
\bd
\lim_{\lai} \frac{1}{l} \log | \E(l,n,2) | = 0 .
\ed
Next, it follows from (\ref{eq:410}) that
\beq
\limsup_{\lai} \sup_{\bzeta \in \GA } \d(l,\bzeta)
& \leq & \sup_{\bzeta \in \GA } - D_c(\bzeta \nm \bmu) \nonumber \\
& = &  - \inf_{\bzeta \in \GA} D_c( \bzeta \nm \bmu ). \nonumber
\eeq
Combining these two facts shows that
\bd
\limsup_{\lai} \frac{1}{l} \log \Pr \{ \bnu(x_1^l) \in \GA \} \leq
- \inf_{\bzeta \in \GA} D_c( \bzeta \nm \bmu ) . \nonumber
\ed
This inequality is in fact stronger than the desired conclusion,
namely the right inequality in (\ref{eq:33}),
because $\GA$ is a subset of $\bar{\GA}$, and as a consequence
\bd
\inf_{\bzeta \in \bar{\GA}} D_c( \bzeta \nm \bmu ) \leq
\inf_{\bzeta \in \GA} D_c( \bzeta \nm \bmu ) ,
\ed
or equivalently
\bd
- \inf_{\bzeta \in \GA} D_c( \bzeta \nm \bmu ) \leq
- \inf_{\bzeta \in \bar{\GA}} D_c( \bzeta \nm \bmu ) .
\ed
Therefore the right inequality in (\ref{eq:33}) is true.
Note that it is pointed out in \cite[pp.\ 16-17]{DZ98} that,
in this special case,
it is possible to replace $\bar{\GA}$ by $\GA$ on the right side
of (\ref{eq:33}).

To establish the left inequality in (\ref{eq:33}),
suppose $\bzeta$ is an interior point of $\GA$.
Now (\ref{eq:47b}) implies that there is a sequence $\{ \bzeta_l \}$ 
where $\bzeta_l \in \E(l,n,2)$ such that $\bzeta_l \ap \bzeta$ as $\lai$.
Moreover, since $\bzeta$ is an interior point of $\GA$,
it can be assumed that $\bzeta_l \in \GA$ for all sufficiently large $l$.
Hence for each sufficiently large integer $l$, we have
\bd
\Pr \{ \bnu(x_1^l) \in \GA \} \geq 
\Pr \{ \bnu(x_1^l) = \bzeta_l \} ,
\ed
\beq
\frac{1}{l} \log \Pr \{ \bnu(x_1^l) \in \GA \} & \geq &
\d(l,\bzeta_l) \nonumber \\
& \geq & - D_c( \bzeta_l \nm \bmu ) + o(1/l) \nonumber \\
& \ap & - D_c( \bzeta \nm \bmu ) \mbox{ as } \kai , \nonumber
\eeq
because $\bzeta_l \ap \bzeta$ as $\lai$ and $D(\cdot \nm \cdot)$ is 
continuous.
Hence it follows that
\bd
\liminf_{\lai} \frac{1}{l} \log 
\Pr \{ \bzeta(x_1^l) \in \GA \} \geq - D_c(\bzeta \nm \bmu ) ,
\fa \bzeta \in \GA^o .
\ed
Since the above inequality holds for {\it every\/} $\bzeta \in \GA^o$, we can
conclude that
\beq
\liminf_{\lai} \frac{1}{l} \log \Pr \{ \bzeta(x_1^l) \in \GA \} & \geq &
\sup_{\bzeta \in \GA^o} - D_c(\bzeta \nm \bmu) \nonumber \\
& = &  - \inf_{\bzeta \in \GA^o} D_c( \bzeta \nm \bmu ) . \nonumber
\eeq
This establishes that the relationships in (\ref{eq:33}) hold with
$I(\bzeta) = D_c( \bzeta \nm \bmu )$.

To complete the proof, it is only necessary to establish that the function
$I(\bzeta) = D_c( \bzeta \nm \bmu )$ is lower-semicontinuous.
But in fact it is continuous.
Hence $D_c( \bzeta \nm \bmu )$ is the rate function.
$\halmos$

In conclusion, it may be remarked that when the samples $\{ x_t \}$
come from an i.i.d.\ process, there are exact formulae for both the size
of each type class as well as the likelihood of each sample within a type class.
In the case where the samples come from a Markov process, there are only
bounds.
{\em However}, the ``correction terms'' in these bounds approach zero
as the number of samples approaches infinity, thus allowing us to
deduce the rate function for doublet frequencies
in a straight-forward fashion.

Next we present the rate fuction for $s$-step Markov processes and the
empirical measure defined in (\ref{eq:61}).
\begin{theorem}\label{thm:61}
Suppose $\{X_t\}$ is a stationary $s$-step Markov assuming values in the
finite set $\Abb$, with the $(s+1)$-tuple frequency vector $\bmu \in \McAs$,
and suppose that $\bmu_\i > 0 \fa \i \in \Abb^{s+1}$.
Define $\bnu(x_1^l) \in \McAs$ as in (\ref{eq:61}).
Then the $\McAs$-valued process $\{ \bnu(x_1^l) \}$ satisfies the LDP
with the rate function
\be\label{eq:63}
I(\bnu) := D_c( \bnu \nm \bmu ) = D( \bnu \nm \bmu ) - D( \bnubar \nm \bmubar ) .
\ee
\end{theorem}

{\bf Proof:}
The proof is the same as that of Theorem \ref{thm:41}, except for more
messy notation (which is why Theorem \ref{thm:41} is stated and proved first).
The first step is to approximate the log likelihood of each sample path
$x_1^l \in \Abb^l$ in terms of the associated empirical distribution
$\bnu(x_1^l)$ defined in (\ref{eq:61}),
and its reduced version $\bnubar(x_1^l)$ defined in the familiar manner,
namely
\bd
\nubar_\i = \sum_{j \in \Abb} \nu_{\i j} , \fa \i \in \Abb^s .
\ed
In analogy with earlier reasoning, using the shorthand
\bd
\Pr \{ x_t | x_{t-s}^{t-1} \} =
\Pr \{ X_t = x_t | X_{t-s}^{t-1} = x_{t-s}^{t-1} \} ,
\ed
we can write
\beq
\Pr \{ X_1^l = x_1^l \} & = & \mubar(x_1^s) \cdot
\prod_{t=s+1}^l \Pr \{ x_t | x_{t-s}^{t-1} \} 
\nonumber \\
& = & \mubar(x_1^s) \cdot
\prod_{t=s+1}^l \frac{ \mu (x^t_{t-s}) }{ \mubar (x^{t-1}_{t-s}) } ,
\nonumber \\
& = & A \cdot B , \nonumber
\eeq
where
\bd
A = \mubar(x_1^s) \cdot
\prod_{t=l+1}^{l+s} \frac{ \mubar ( x_{t-s}^{t-1} ) }{ \mu ( x_{s-1}^t ) } ,
\ed
\bd
B = \prod_{t=s+1}^{l+s} \frac { \mu ( x_{s-1}^t ) } { \mubar ( x_{t-s}^{t-1} ) } .
\ed
It is easy to verify that $A$ is just the generalization of
$(\mubar(x_1) \mubar(x_l)) / \mu(x_1 x_l)$ to $s$-step Markov processes.
Now, as in the proof of Theorem \ref{thm:41}, it is possible to bound
the term $A$ both above and below.
Since $\mu_\i > 0 $ for all $\i \in \Abb^{s+1}$, choose constants
$\underline{a}, \bar{a}, \underline{b}, \bar{b}$ such that
\bd
0 < \underline{a} \leq \mubar_\i \leq \bar{a} \fa \i \in \Abb^s ,
\ed
\bd
0 < \underline{b} \leq \mu_\i \leq \bar{b} \fa \i \in \Abb^{s+1} .
\ed
Now let us define constants $\underline{c},\bar{c}$ by
\bd
\underline{c} = (s+1) \log \underline{a} - s \log \bar{b} ,
\ed
\bd
\bar{c} = (s+1) \log \bar{a} - s \log \underline{b} .
\ed
Then
\bd
\underline{c} \leq \log A \leq \bar{c} .
\ed

Next let us compute $\log B$.
This leads to
\bd
\log B = \sum_{t=s+1}^{l+s} [ \log \mu ( x_{s-1}^t ) - \log \mubar ( x_{t-s}^{t-1} ) ] .
\ed
Recall that we had earlier defined $l_\i = l \nu_\i$ for all $\i \in \Abb^{s+1}$
and $\lbar_\i = l \nubar_\i$ for all $\i \in \Abb^s$.
So in the above summation, the string $x_{s-1}^t$ equals $\i \in \Abb^{s+1}$
exactly $l_\i$ times, while the string $x_{t-s}^{t-1}$ equals $\i \in \Abb^s$
exactly $\lbar_\i$ times.
Therefore
\beq
\log B & = & l \left[ \sum_{\i \in \Abb^{s+1} } \nu_\i \log \mu_\i 
- \sum_{\i \in \Abb^s} \nubar_\i \log \mubar_\i \right] \nonumber \\
& = & - l [ J( \bnu,\bmu) - J(\bnubar,\bmubar) ] . \nonumber
\eeq
Therefore (\ref{eq:49a}) and (\ref{eq:49b}) hold with the updated
definitions of the constants $\underline{c}$ and $\bar{c}$.
The rest of the proof exactly follows that of Theorem \ref{thm:41}.
$\halmos$

\section{The Rate Function for Singleton Frequencies}\label{sec:6}

In this section we use Theorems \ref{thm:41} and \ref{thm:61},
together with a very general result known as the ``contraction principle,''
to derive the rate function for singleton frequencies of (one-step)
Markov chains.
This is the approach adopted in \cite[Chapter IV]{Hollander00}.
In some earlier work such as \cite{DV75b,DV75c,DV76}, 
the rate function for singleton frequencies of a Markov processes 
(Theorem \ref{thm:32}) is derived first, and this is used
to derive the rate function for doublet frequencies (Theorem \ref{thm:33}).
This is also the approach followed in \cite{DZ98}.

Now we state a special case of the contraction principle, as needed
for our work.
The reader is referred to \cite[Theorem III.20]{Hollander00} or
\cite[Theorem 4.2.1]{DZ98} for a general statement of the result.
In particular, the general result requires that the rate function
function $I(\cdot)$ should, in addition to being lower semi-continuous,
also be a ``good'' rate function in the sense that it has compact
level sets.
In other words, for every constant $c$, the set $\{ \bnu : I( \bnu ) \leq c \}$
should be a compact set.
In the special situation being studied here, where the stochastic
processes evolve over a finite alphabet, the set of measures is a subset 
of a simplex $\Sm_m$ of an appropriate dimension, which is itself a
compact set.
Thus {\it all\/} rate functions are ``good.''
But in case the reader wishes to study more general situations,
this additional requirement, which is automatically satisfied in
the present setting, must be taken into account.

With that caveat, we state the contraction principle as it
applies to the current setting.
Consider first the case of a one-step Markov process.
Define a map $\f: \McA \ap \Sm_n$ by
\bd
[ \f (\bnu) ]_i := \sum_{j \in \Abb} \nu_{ij} = \sum_{j \in \Abb} \nu_{ji} .
\ed
Thus $\f$ maps a stationary distribution $\bnu$ on $\Abb^2$ onto its
one-dimensional marginal $\bnubar$.
Moreover, if we construct $\bnu(x_1^l)$ for a sample $x_1^l$ using
the formula (\ref{eq:41}), then the corresponding $\f[\bnu(x_1^l)]$ is
the usual empirical distribution of singleton frequencies $\bphi$
defined in (\ref{eq:34a}).
Now, by invoking the contraction principle (see e.g.\ \cite{DZ98,FK06}),
we can readily conclude the following:

\begin{theorem}\label{thm:51}
The $\Sm_n$-valued process $\bphi(x_1^l)$ satisfies the LDP with the
rate function
\be\label{eq:51}
J(\bphi) := \inf_{\bnu \in \McA} D_c( \bnu \nm \bmu ) \st \bnubar = \bphi .
\ee
\end{theorem}

Recall that
\bd
D_c( \bnu \nm \bmu ) = D( \bnu \nm \bmu ) - D( \bnubar \nm \bmubar ) .
\ed
Hence we can also write
\be\label{eq:52}
J(\bphi) = \left[ \inf_{\bnu \in \McA} D( \bnu \nm \bmu ) \st 
\bnubar = \bphi \right] - D( \bphi \nm \bmubar ) .
\ee

The remainder of the proof follows \cite[pp.\ 45-47]{Hollander00}.
For the convenience of the reader, the relevant material is
reproduced here.

\begin{theorem}\label{thm:51a}
The quantity $J(\bphi)$ defined in (\ref{eq:52}) is also given by
\be\label{eq:53}
J(\bphi) = \sup_{\u > \bz} \sum_{i=1}^n \phi_i \log \frac{u_i}{( A \u )_i} .
\ee
\end{theorem}

{\bf Proof:}
Let us define
\bd
g(\u) := \sum_{i=1}^n \phi_i \log \frac{u_i}{( A \u )_i} .
\ed
Then the claim is that the supremum of $g(\cdot)$ equals the infimum
in (\ref{eq:51}).
To show this, observe that $g(\u) \ap \infty$ if any component of $\u$
approaches zero.
Thus the supremum is achieved at some point in the interior of the positive
orthant.
The conditions characterizing the supremum can be obtained using elementary
calculus.
We have
\bd
\frac{ \partial g}{ \partial u_i } = 
\frac{ \phi_i }{ u_i } - \sum_{j \in \Abb} \frac{ \phi_j a_{ji} }{( A \u )_j} .
\ed
Thus any optimum $\u^*$ satisfies\footnote{Clearly the optimum cannot
be unique because any multiple of $\u^*$ is also optimal.}
\be\label{eq:54}
\phi_i = \sum_{j \in \Abb} \frac{ \phi_j a_{ji} u^*_i }{( A \u^* )_j} .
\ee
Define
\be\label{eq:55}
b_{ji} = \frac{ a_{ji} u^*_i }{( A \u^* )_j} , \mbox{ or }
b_{ij} = \frac{ a_{ij} u^*_j }{( A \u^* )_i} ,
\ee
and note that the matrix $B = [b_{ij} ]$ is stochastic in that
$\sum_{j} b_{ij} = 1$ for all $i$.
Now (\ref{eq:54}) can be interpreted as $\bphi = \bphi B$.
So if we define $\nu^*_{ij} := \phi_i b_{ij}$, it follows that
$\bnu^* \in \McA$ and that $\bar{\bnu^*} = \bphi$.

Next, let us compute $D_c(\bnu^* \nm \bmu)$.
From the earlier proof that the quantities in (\ref{eq:43a}) and (\ref{eq:47})
are equal, we know that
\beq
D_c(\bnu^* \nm \bmu)
& = & \sum_{i} \phi_i \sum_j b_{ij} \log \frac{ b_{ij} }{ a_{ij} } \nonumber \\
& = & \sum_{i} \sum_j \nu^*_{ij} \log \frac{ b_{ij} }{ a_{ij} } . \nonumber
\eeq
Now (\ref{eq:55}) shows that $b_{ij}/a_{ij} = u^*_j/(A \u^*)_i$.
Therefore
\bd
D_c(\bnu^* \nm \bmu) = \sum_{i} \sum_j \nu^*_{ij} [ \log u^*_j 
- \log (A \u^* )_i ] .
\ed
Now we prove a handy identity:
Let $c_1 , \ldots , c_n$ be any constants.
Then
\be\label{eq:56}
\sum_i \sum_j \nu^*_{ij} c_i = \sum_{i} \sum_j \nu^*_{ij} c_j .
\ee
To prove this, note that
\bd
{\rm LHS } = \sum_i \phi_i c_i = \sum_j \phi_j c_j = {\rm RHS}.
\ed
Therefore 
\beq
D_c(\bnu^* \nm \bmu) & = & \sum_{i} \sum_j \nu^*_{ij} [ \log u^*_i
- \log (A \u^* )_i ] \nonumber \\
& = & \sum_{i} \sum_j \nu^*_{ij} \log \frac{ u^*_i }{ (A \u^* )_i } 
\nonumber \\
& = & \sum_i \phi_i \log \frac{ u^*_i }{ (A \u^* )_i } = g(\u^*) . \nonumber
\eeq

Now suppose that $\bnu \in \McA$ satisfies $\bnubar = \bphi$.
It is now shown that
\bd
D_c( \bnu \nm \bmu ) \geq g(\u^*) = D_c(\bnu^* \nm \bmu) .
\ed
Once the above relationship is established, it then follows
from the contraction principle that $D_c(\bnu^* \nm \bmu) = g(\u^*)$
is indeed the minimum in (\ref{eq:51}) and is thus the rate function
$J(\bphi)$.
For this purpose, we make use of (\ref{eq:47d}) to compute the quantity
$D_c(\bnu^* \nm \bmu)$.
As before, let $c_{ij} = \nu_{ij} / \nubar_i$, and 
$a_{ij} = \mu_{ij} / \mubar_i$.
Also note from (\ref{eq:55}) that $b_{ij} / a_{ij} = u^*_j / (A \u^*)_i $.
With these observations, we reason as follows:
\beq
D_c( \bnu \nm \bmu ) & = & \sum_i \sum_j \nu_{ij} \log \frac{c_{ij}}{a_{ij}}
\nonumber \\
& = & \sum_i \sum_j \nu_{ij} \log \frac{c_{ij}}{b_{ij}} \nonumber \\
& + & \sum_i \sum_j \nu_{ij} \log \frac{b_{ij}}{a_{ij}} \nonumber \\
& = & D_c( \bnu \nm \bnu^* ) \nonumber \\
& + & \sum_{i} \sum_j \nu_{ij} [ \log u^*_j - \log (A \u^* )_i ] \nonumber \\
& = & D_c( \bnu \nm \bnu^* ) \nonumber \\
& + & \sum_{i} \phi_i [ \log u^*_i - \log (A \u^* )_i ] \nonumber \\
& = & D_c( \bnu \nm \bnu^* ) + g(\u^*) , \nonumber 
\eeq
where in the next to last step we interchange $u_j$ with $u_i$ using
(\ref{eq:56}) and use the fact that $\bnubar = \bphi$.
Since $D_c( \bnu \nm \bnu^* ) \geq 0$, it follows that
\bd
D_c( \bnu \nm \bmu ) \geq g(\u^*) \mbox{ if } \bnubar = \bphi .
\ed
Thus the quantity $J(\bphi)$ defined in (\ref{eq:51}) is also given by
the supremum in (\ref{eq:53}).
$\halmos$

Theorem \ref{thm:61} presents the rate function for $(s+1)$-tuple
frequencies of an $s$-step Markov chain.
Using the contraction principle, it is possible to obtain the rate function
for the frequencies of $k$-tuples where $1 \leq k \leq s$.
However, perhaps that is an example of ``the principle of diminishing returns,''
as the notation gets truly complex without a great deal of fresh insight.
So this generalization is left to the reader.

\section{Conclusions}\label{sec:9}

In this paper, we have derived several known results on the large deviation
propery (LDP) of Markov chains using only elementary arguments based
on the method of types, and avoiding advanced arguments based on the
G\"{a}rtner-Ellis theorem, as in standard texts such as \cite{DZ98},
or Varadhan's lemma as in \cite{Hollander00}.
As a result, the proofs are far more accessible.
Moreover, the extension to multi-step Markov chains
is quite transparent and easy.
The various rate functions are given a natural interpretation in terms
of differential relative entropy.

\section*{Appendix A: Upper and Lower Bounds for the Cardinality
of Type Classes}

In this appendix we state and prove upper and lower bounds for the
cardinality of type classes.
Recall that if $x_1^l$ is a sample path in $\Abb^l$, then for
one-step Markov chains the associated empirical distribution
is defined by (\ref{eq:41}), while the associated empirical distribution
for $s$-step Markov chains is defined by (\ref{eq:61}).
The symbol $\E(l,n,s+1)$ denotes the set of empirical distributions
in $\McAs$ that correspond to sample paths in $\Abb^l$, and for
$\bzeta \in \E(l,n,s+1)$, the symbol $T(\bzeta,l,s+1)$ denotes
the set of all sample paths $x_1^l \in \Abb^l$ such that the 
associated empirical distribution $\bnu(x_1^l)$ equals $\bzeta$,
and is known as the type class of $\bzeta$.
The objective here is to derive upper and lower bounds for the cardinality
$| T(\bzeta,l,s+1) |$ in terms of the properties of $\bzeta$.
Slightly less tight bounds for sufficiently long sample paths
are also quoted in \cite[(5.5) and (5.6)]{CG05}.
Note that these bounds are based on a line of reasoning that goes
back to \cite{DLS81}; however, the present bounds are less complicated.
Essentially the same bounds are derived in \cite[p.\ 17]{Hollander00} in
an extremely terse form.
Since the discussion in \cite{Hollander00} does not fully cover all the
intricacies, we present the proof in some detail.

\begin{theorem}\label{thm:A1}
Suppose $\bzeta \in \E(l,n,s+1)$.
Then the cardinality of the type class $T(\bzeta,l,2)$ is bounded by
\be\label{eq:A1}
(2l)^{-n^{s+1}} e^{lD(\bzeta)} \leq | T(\bzeta,l,s+1)| \leq l e^{lD(\bzeta)}
\fa l \geq n .
\ee
\end{theorem}

{\bf Proof:}
In order to prove the theorem, we begin by introducing a symbol for
the combinatorial index when there are $n$ possible outcomes.
Suppose $m_1 , \ldots , m_n$ are integers with $m := \sum_{i=1}^n m_i$.
For notational convenience let $\m$ denote the vector $[m_i]$.
Then the integer
\bd
C^m_{\m} := \frac{ m! }{\prod_{i=1}^n m_i ! }
\ed
represents the number of distinct ways of assigning $m$ labels to the $n$
elements of $\Abb$ in such a way that the $i$-th element is assigned
precisely $m_i$ labels.
If $n=2$ so that $m_2 = m - m_1$, then it is customary to write just
$C^m_{m_1}$ or $C^m_{m_2}$.
Define the associated probability distribution
\bd
\bzeta := [ m_i/m , i = 1 , \ldots , n ] \in \Sm_n .
\ed
Then a classic result found for example in \cite[Lemmas 1 and 2]{Csiszar98}
states that
\be\label{eq:A2}
\left( C_{n-1}^{m+n-1} \right)^{-1} \exp[ H(\bzeta) ]
\leq C^m_{\m} \leq \exp[ H(\bzeta) ] .
\ee
Since we are interested in the case where $m$ approaches infinity, we
observe that, for all $m \geq n$, we have
\beq
C_{n-1}^{m+n-1} & = & \frac{ (m - n - 1)!}{ (n-1)! m!} =
\frac { \prod_{i=1}^{n-1} ( m + n - i ) } { (n-1)! } \nonumber \\
& \leq & (2m)^{n-1} . \nonumber
\eeq
As a result we can rewrite (\ref{eq:A2}) as
\be\label{eq:A3}
(2m)^{-(n-1)} \exp[ H(\bzeta) ] \leq C^m_{\m} \leq \exp[ H(\bzeta) ] ,
\fa m \geq n .
\ee
This bound is less sharp than that in (\ref{eq:A2}), but is easier to
work with.
As we shall see, so long as we get the exponential term exactly right,
the polynomial terms in front don't really matter.

To prove the estimate (\ref{eq:A1}), we begin with the case $s = 1$ so that
we are estimating doublet frequencies, and then generalize.
This will result in less messy notation at the outset, and lead to greater
insight.

Suppose $\bzeta \in \E(l,n,2)$.
As we have done elsewhere, we define $\bzetabar \in \Sm_n$ as the reduced
version of $\bzeta$, and $l_{ij} = l \zeta_{ij} , \lbar_i = l \zetabar_i$
to be associated integers.
With the distribution $\bzeta$, we can associate a
directed graph $\G(\bzeta)$ with $n$ nodes and $l$ edges by placing
$l_{ij}$ edges from node $i$ to node $j$.
Because $\bzeta$ is a stationary distribution, it is clear that every
node in the graph $\G(\bzeta)$ has equal in-degree and out-degree.
Therefore the graph $\G(\bzeta)$ is a union of cycles.
In other words, it is possible to start at a node, and trace out a path
that ends at the same node.
Note that it is permissible for the path
to pass through any of the $n$ nodes one or more times.
Indeed if $l > n$, this is unavoidable.
The cardinality of the associated type class $T(\bzeta,l,2)$ is
precisely equal to the number of different ways of tracing out such a path.
While every $\bzeta \in \McA$ corresponds to a graph where every node
has equal in- and out-degrees, not every $\bzeta \in \McA$ (or equivalently,
not every graph with this property) belongs $\E(l,n,2)$.
For example, the distribution 
\bd
\bzeta_{ij}  = \frac{ \lbar_i }{ l } \d_{ij} , \fa i, j \in \Abb ,
\ed
with $\lbar_i > 0$ for all $i$, where $\d$ is the Kronecker delta,
belongs to $\McA$ but not to $\E(l,n,2)$.
In order for $\bzeta$ to belong to $\E(l,n,2)$ and not just
$\McA$, a necessary and sufficient condition is that the graph 
$\G(\bzeta)$ should consist of one strongly connected component.
This explains why the ``diagonal'' distribution above fails to belong to
$\E(l,n,2)$, because in this case there
is more than one strongly connected component.

Next we introduce a way of describing a graph $\G(\bzeta)$ associated
with a $\bzeta \in \McA$.

Recall that $\{ a_1 , \ldots , a_n \}$ denote the elements of the state space
$\Abb$, ordered in some arbitrary fashion.
Given a sample path $x_1^l$ and the associated empirical estimate
$\bnu(x_1^l)$ constructed as in (\ref{eq:41}), let us define $n$
sequences $S(1)$ through $S(n)$ as follows:
The set $S(i)$ has cardinality $\lbar_i$, and consists of all the symbols
that follow $a_i$ in the sample path, in that order.
To illustrate, suppose $n = 3$, $\Abb = \{ a, b, c \}$\footnote{This is
clearer than writing $A = \{ 1, 2, 3 \}$ or
$\Abb = \{ a_1 , a_2 , a_3 \}$.}, $l = 10$, and 
\bd
x_1^l = abaccbacbc .
\ed
Then, with the ghost transition from $c$ back to $a$ added, the three
sets are defined as
\bd
S(a) = bcc , S(b) = aac, S(c) = cbba .
\ed
It is easy to see that, given the $n$ sets $S(1)$ through $S(n)$,
we can reconstruct the associated cycle;
however, we would not know the starting point.
In other words, given the sets $S(1)$ through $S(n)$, there are at most
$n$ different corresponding sample paths, corresponding to choosing one
of the $n$ nodes as the starting point of the cycle.
Moreover, if $\bth$ belongs to the same type class as $\bzeta$,
then each $S(i,\bth)$ is a permutation of the corresponding set $S(i,\bzeta)$.

This suggests a way of finding an upper bound for $ | T(\bzeta,l,2) |$.
Given a sample path $x_1^l$ and the associated distribution $\bzeta$,
in order to enumerate all elements of $T(\bzeta,l,2)$,
there are two things we can do.
First, we can permute the elements of $S(i)$ for $i = 1 , \ldots , n$,
which leads to $\prod_{i=1}^n S(i)!$ variants of the sample path that
generate the same empirical distribution.
Second, for each such permutation we can do a cyclic shift of the
sample path and thus change the starting point.
Note all permutations of each $S(i)$ and/or not all cyclic shifts of
the starting point lead to a valid sample path, or to distinct sample
paths.
But the combination of these two actions provides {\it an upper bound\/}
on the cardinality of $ | T(\bzeta,l,2) |$.
Specifically, it follows that
\beq
|T(\bzeta,l,2)|
& \leq & n \frac{\prod_{i \in \Abb} \lbar_i !}
{\prod_{i \in \Abb} \prod_{j \in \Abb} l_{ij} ! } \nonumber \\
& \leq & l \frac{\prod_{i \in \Abb} \lbar_i !}
{\prod_{i \in \Abb} \prod_{j \in \Abb} l_{ij} ! }  \fa l \geq n . \label{eq:A3a}
\eeq
The denominator term arises because each of the $l_{ij}$ edges
are indistinguishable for each $i,j$, so we need to divide by
the number of ways in which these can be permuted.

Finding a lower bound for $ | T(\bzeta,l,2) |$ is a bit more involved.
To obtain a lower bound on the number of paths,
we use the following argument from \cite{DLS81},
which is explained more clearly in \cite[p.\ 17]{Hollander00}.
Pick any one cycle that spans each node in the connected component exactly once,
and then delete all of these edges from the graph.
Delete the corresponding elements from the sets $S(1)$ through $S(n)$.
This has the effect of reducing the cardinality of each set $S(i)$ by
exactly one.
Then all possible permutations of these reduced sets will result in
valid paths, because it is not possible to get `stuck' at any node --
the edges of the deleted cycle serve as an escape route.
Thus there at least $\prod_{i \in \Abb} (\lbar_i - 1)!$ permutations
that will result in paths.
Therefore
\be\label{eq:A3b}
\frac{\prod_{i \in \Abb} ( \lbar_i - 1) !}
{\prod_{i \in \Abb} \prod_{j \in \Abb} l_{ij} ! }
\leq |T(\bzeta,l,2)|
\ee
Again, the denominator term arises because each of the $l_{ij}$ edges
are indistinguishable for each $i,j$, so we need to divide by
the number of ways in which these can be permuted.
Combining (\ref{eq:A3a}) and (\ref{eq:A3b}) leads to the desired
conclusion (\ref{eq:A1}) when $s = 1$.

For multi-step Markov chains, suppose $\bzeta \in \E(l,n,s+1)$.
Then $l \bzeta$ has only integer components, as does $l \bzetabar$.
we can associate a directed graph with $\bzeta \in \McAs$ as follows:
The graph has $n^s$ nodes labelled as the elements of $\Abb^s$.
Suppose $i, k \in \Abb, \j \in \Abb^{s-1}$.
Then between nodes $i \j$ and $\j k$, we draw $l \zeta_{i \j k}$ directed edges.
Then the in-degree of node $\j k$ equals $\sum_{i \in \Abb}
l \zeta_{i \j k}$, which we denote as before by $\lbar_{\j k}$.
With this argument, it follows that for each $\bzeta \in \E(l,n,s+1)$, we have
\beq
\frac { \prod_{\i \in \Abb^s} ( \lbar_{\i} - 1)! }
{ \prod_{\i \in \Abb^s} \prod_{j \in \Abb} l_{\i j}! }
& \leq & | T(l,n,s+1) | \nonumber \\
& \leq & l \frac { \prod_{\i \in \Abb^s} \lbar_{\i}! }
{ \prod_{\i \in \Abb^s} \prod_{j \in \Abb} l_{\i j}! } . \label{eq:A4}
\eeq
In the above, if $l_{\i j} = 0$ for some $\i, j$, then that term is
omitted from the product, and of course we take $0! = 1$.
Now observe that
\bd
( \lbar_{\i} - 1)! = \frac{ \lbar_\i ! }{ \lbar_\i } 
\geq \frac { \lbar_\i ! }{ l } , \fa \i \in \Abb^s .
\ed
Hence we can rewrite (\ref{eq:A4}) as
\beq
\frac{1}{l^{n^s}} \frac { \prod_{\i \in \Abb^s} \lbar_{\i}! }
{ \prod_{\i \in \Abb^s} \prod_{j \in \Abb} l_{\i j}! }
& \leq & | T(l,n,s+1) \nonumber \\
& \leq & l \frac { \prod_{\i \in \Abb^s} \lbar_{\i}! } { \prod_{\i \in \Abb^s}
\prod_{j \in \Abb} l_{\i j}! } . \label{eq:A5}
\eeq
Since $\sum_{j \in \Abb} l_{\i j} = \lbar_\i$ for all $\i \in \Abb^s$,
we can see lots of occurrances of the combinatorial parameter in (\ref{eq:A4}).
We can rewrite (\ref{eq:A4}) as
\be\label{eq:A6}
\frac{1}{l^{n^s}} \prod_{\i \in \Abb^s} C^{\lbar_\i}_{l_{\i 1 , \ldots , \i n}}
\leq | T(l,n,s+1) | \leq
l \prod_{\i \in \Abb^s} C^{\lbar_\i}_{l_{\i 1 , \ldots , \i n}} .
\ee
Now we make use of the upper and lower bounds in (\ref{eq:A4}).
For this purpose, let us define
\bd
\bzetabar_\i := [ l_{\i 1} / \lbar_\i \dots l_{\i n} / \lbar_\i ] \in \Sm_n 
\fa \in \Abb^s .
\ed
Then it follows from (\ref{eq:A4}) that, for all $\i \in \Abb^s $,
\beq
\frac{1}{2^{n-1} \lbar_\i^{n-1} } \exp ( \lbar_\i H( \bzetabar_\i ) )
& \leq & C^{\lbar_\i}_{l_{\i 1 , \ldots , \i n}} \nonumber \\
& \leq & \exp ( \lbar_\i H( \bzetabar_\i ) ) . \label{eq:A7}
\eeq
When we substitute these bounds in (\ref{eq:A6}), we need compute two
quantities, namely
\bd
l^{n^s} \prod_{\i \in \Abb^s} 2^{n-1} \lbar_\i^{n-1} , 
\ed
and
\bd
\prod_{\i \in \Abb^s} \exp ( \lbar_\i H( \bzetabar_\i ) )
= \exp \left( \sum_{\i \in \Abb^s} \lbar_\i H( \bzetabar_\i ) \right) .
\ed
The first term is easy.
Since $\sum_{\i \in \Abb^s} \lbar_\i = l^s$, the first product is
\bd
(2^{n-1} l^n)^{n^s} \leq (2l)^{n^{s+1}} .
\ed
As for the second term, we get
\beq
\lbar_\i H( \bzetabar_\i ) & = & - \lbar_\i 
\sum_{j \in \Abb} \frac{ l_{\i j} }{ \lbar_\i }
[ \log l_{\i j} - \log \lbar_\i ] \nonumber \\
& = & \lbar_\i \log \lbar_\i - \sum_{j \in \Abb} l_{\i j} \log l_{\i j} ,
\nonumber
\eeq
where we use the fact that $\sum_{j \in \Abb} l_{\i j} = \lbar_\i$.
Therefore
\beq
\sum_{\i \in \Abb^s} \lbar_\i H( \bzetabar_\i ) & = &
- \sum_{\i \in \Abb^s} \sum_{j \in \Abb} l_{\i j} \log l_{\i j} \nonumber \\
& + & \sum_{\i \in \Abb^s} \lbar_\i \log \lbar_\i \nonumber \\
& = & l [ H( \bzeta ) - H( \bzetabar) ] = l H_c( \bzeta ) . \label{eq:A8}
\eeq
Substituting from (\ref{eq:A7}) and (\ref{eq:A8}) into (\ref{eq:A6})
leads to the desired bound, namely
\bd
(2l)^{-n^{s+1}} \exp ( l D(\bzeta) ) \leq | T(\bzeta,l,s+1) | \leq
l \exp ( l D(\bzeta) ) .
\ed

\end{document}